\documentclass[12pt]{amsart}
\usepackage{amsmath,amssymb,amscd,amsthm,amsxtra}
\newtheorem{Thm}{Theorem}[section]

\newtheorem{Cor}[Thm]{Corollary}
\newtheorem{Lem}[Thm]{Lemma}
\newtheorem{Prop}[Thm]{Proposition}

\def\ldots{\mathinner{\ldotp\ldotp\ldotp}}
\def\ldots{\mathinner{\cdotp\cdotp\cdotp}}

\def \Bbb{\mathbb}

\def\Ave{\mathop{\text {Ave}}}
\def\supp{\text{supp }}

\begin{document}

\title{The Marcinkiewicz multiplier condition for
bilinear operators}

\author{Loukas Grafakos}
\address{Department of Mathematics \\
University of Missouri-Columbia \\
Columbia, MO 65211
}
\email{loukas@math.missouri.edu}

\author{Nigel J. Kalton}
\address{Department of Mathematics \\
University of Missouri-Columbia \\
Columbia, MO 65211
}
\email{nigel@math.missouri.edu}

\thanks{The research of both authors was supported by the NSF}

\date{\today}



\begin{abstract}
This article is concerned with the question of whether
Marcin- kiewicz multipliers on $\mathbb R^{2n}$ give
rise to   bilinear multipliers
on $\mathbb R^n\times \mathbb R^n$. We show that  this is not always the
case.  Moreover we find necessary and sufficient  conditions for such bilinear
multipliers to be bounded.  These conditions in particular imply that a slight
logarithmic   modification of the Marcinkiewicz condition  gives
multipliers for which the corresponding  bilinear operators are bounded on
products of Lebesgue   and   Hardy spaces.
\end{abstract}
\maketitle

\section{Introduction}\label{s-introduction}
 \setcounter{equation}0

In this article we study bilinear multipliers of Marcinkiewicz type.
Recall that a  function $\sigma(\xi, \eta)= \sigma(\xi_1, \dots , \xi_n, \eta_1,
\dots , \eta_n)$ defined    away from the coordinate axes  on $\mathbb R^{2n}$,
 which satisfies  the conditions
\begin{equation}\label{ma}
|\partial^\alpha_\xi \partial^\beta_\eta \sigma(\xi, \eta ) | \le
C_{\alpha, \beta} |\xi_1|^{- \alpha_1}\dots |\xi_n|^{- \alpha_n}
|\eta_1|^{- \beta_1}\dots |\eta_n|^{- \beta_n}
\end{equation}
for sufficiently large multi-indices $\alpha=(\alpha_1, \dots , \alpha_n)$ and
$\beta=(\beta_1, \dots , \beta_n)$, is called a Marcinkiewicz multiplier.
It is a classical result, see for instance \cite{stein-old}, that Marcinkiewicz
multipliers give  rise to   bounded linear
operators $M_\sigma$ from  $L_p(\mathbb R^{2n})$ into itself for $1<p<\infty$.
Here $M_\sigma$ is the multiplier operator with symbol   $\sigma$, that is
$$
M_\sigma (F)(x) = \int_{\mathbb R^{2n}} \widehat{F}(\xi) \sigma (\xi) e^{ 2\pi i
\langle x, \xi \rangle } d\xi,
$$
where $F$ is a Schwartz function on $\mathbb R^{2n}$ and
$\widehat{F}(\xi)$ is the
Fourier transform of $F$, defined by
$
\widehat{F}(\xi)= \int_{\mathbb R^{2n}} F(x)  e^{ -2\pi i
\langle x, \xi \rangle } dx.
$
(We will use the notation
$\langle x, y \rangle=\sum_{k=1}^m  x_ky_k$ for $x=(x_1, \dots , x_m)$ and
$y=(y_1, \dots , y_m)$ elements of $\mathbb R^{m}$.)
The Marcinkiewicz  condition (\ref{ma}) is less restrictive  than
the H\"ormander-Mihlin condition
\begin{equation}\label{ho}
|\partial^\alpha_\xi \partial^\beta_\eta \sigma(\xi, \eta ) | \le
C_{\alpha, \beta} (|\xi|+|\eta|)^{-|\alpha|  -|\beta|},
\end{equation}
which  is also known to imply boundedness for the linear
operator $W_\sigma$ from  $L_p(\mathbb R^{2n})$ into itself when $1<p<\infty$.
The advantage of   condition (\ref{ho}) is that it is supposed to hold
for multi-indices up to   order $|\alpha|+|\beta|\le n+1$ versus
up to order $|\alpha|+|\beta|\le 2n$ for condition  (\ref{ma}). 

In this paper we study bilinear multiplier operators whose symbols
satisfy similar conditions. More precisely,
we are interested in boundedness properties of bilinear operators
$$
W_\sigma (f,g) (x) = \int_{\mathbb R^{2n}} \widehat{f}(\xi)
 \widehat{g}(\eta) \sigma (\xi, \eta ) e^{ 2\pi i
\langle x, \xi  \rangle }e^{ 2\pi i
\langle x,  \eta\rangle } \, d\xi\, d\eta ,
$$
originally defined for
$f,g$   Schwartz functions on $\mathbb R^n$ and
 $\sigma$   a   function on $\mathbb R^{2n}$.   A well-known theorem
of Coifman and Meyer \cite{CM} says that if the function $\sigma$ on
$\mathbb R^{2n}$ satisfies (\ref{ho}) for  sufficiently large multi-indices $\alpha$
and $\beta$, then the bilinear map $W_\sigma (f,g)$  extends to
a bounded operator from $L_{p_1}(\mathbb R^n)\times L_{p_2}(\mathbb R^n)$ into
$L_{p_0,\infty }( \mathbb R^n)$ when $1<p_1, p_2<\infty$, $1/p_1+1/p_2=1/p_0$ and
$p_0\ge 1$.
($L_{p_0,\infty}$  here denotes the space
weak $L_{p_0}$.)
 This result was later extended to the range $1> p_0 \ge 1/2$  by
Grafakos and Torres  \cite{GT} and Kenig and Stein \cite{KS}.
The extension into $L_{p_0}$ for $p_0< 1$ was stimulated by the recent work
of  Lacey and Thiele   \cite{lacey-thiele2} who showed that
the discontinuous symbol $\sigma(\xi, \eta)=-i \text{sgn }(\xi-\eta) $  on
$\mathbb R^2$ gives rise to a bounded bilinear operator $W_\sigma$ from
$L_{p_1}(\mathbb R)\times L_{p_2}(\mathbb R)$ into $ L_{p_0}(\mathbb R)$  for
$2/3<p_0<\infty$ when
$1<p_1,p_2<\infty$ and $1/p_1+1/p_2=1/p_0.$

In this article we address the question of whether  the
Marcinkiewicz condition (\ref{ma}) on  $\mathbb R^{2n}$ gives rise to a  bounded
bilinear operator $W_\sigma$ on  $\mathbb R^n \times \mathbb R^n$.
We answer this question negatively.  More precisely,  we show that there
exist examples of bounded functions $\sigma (\xi, \eta)$ on
$\mathbb R^n \times \mathbb R^n$ which satisfy  the stronger condition
\begin{equation} \label{hypothesis}
|\partial^\alpha_\xi \partial^\beta_\eta \sigma(\xi, \eta ) | \le
C_{\alpha, \beta} |\xi|^{-|\alpha|}|\eta|^{-|\beta|}
\end{equation}
for all multi-indices $\alpha$ and $\beta$,
for which the corresponding bilinear operators $W_\sigma$ do not map
$L_{p_1} \times L_{p_2} $ into
 $L_{p_0,\infty}$ for any triple  of exponents satisfying $1/p_1+1/p_2=1/p_0$
and $1<p_1, p_2<\infty$.

We reduce   this problem to the study of bilinear
operators of the type
\begin{equation}\label{ooo}
(f,g)\to \sum_{j\in \mathbb Z}\sum_{k\in \mathbb Z} a_{jk}\,
\widetilde{\Delta}_jf \, \widetilde{\Delta }_kg  ,
\end{equation}
where $a_{jk}$ is a bounded sequence of scalars depending on $\sigma$ and
$\widetilde{\Delta}_j$ are  the Littlewood-Paley operators given by multiplication
on the Fourier  transform side by a smooth bump supported near the frequency
$|\xi|\sim 2^j$.
In section \ref{s-bilinear+infmatrices},
in particular Theorem \ref{matrixnorms}, we find a necessary and
sufficient  condition on the infinite matrix $A=(a_{jk})_{j,k}$
so that the bilinear
operator in (\ref{ooo}) maps $L_{p_1} \times L_{p_2} $ into
$L_{p_0,\infty}$.
This  condition  is expressed in terms of an Orlicz
space norm of the sequence $(a_{jk})_{j,k}$.
 It turns out that this condition is independent
of  the exponents $p_1,p_2,p_0$ and depends only on quantities intrinsic
to the matrix $A$, (although the actual norm of the operator in (\ref{ooo})
from $L_{p_1} \times L_{p_2} $ into
$L_{p_0,\infty}$  does  depend  on the  indices $p_1,p_2,p_0$).

The results of section \ref{s-bilinear+infmatrices} are transferred to multiplier
theorems for bilinear operators in section \ref{s-applications}.
This transference is achieved using a Fourier expansion of the symbol $\sigma$
on products of dyadic cubes.
Theorem \ref{main1} is the main result of this section and Theorem
\ref{bestposs} shows that this theorem is best possible.
Theorem  \ref{main1} allows us to derive  that  the estimates
\begin{equation}\label{above}
 |\partial_{\xi}^{\alpha}\partial_{\eta}^{\beta} \sigma(\xi,\eta)| \le
C_{\alpha, \beta }|\xi|^{-|\alpha|} |\eta|^{-|\beta|}
\big(\log(1+|\log\tfrac{|\xi|}{|\eta|}|)\big)^{-\theta}
\end{equation}
do  give  rise to a bounded bilinear operator $W_\sigma$  on
products of $L_p$ spaces when $\theta >1$, while we   show that this is not
the case when $0<\theta<\frac12.$  We obtain similar results when
the expression $\big(\log(1+|\log\tfrac{|\xi|}{|\eta|}|)\big)^{-\theta}$
in (\ref{above})  is
replaced by the  expression $\big(\log(1+|\log\tfrac{|\xi|}{|\eta|}|)\big)^{-1}
 \!\big(\log(1+\log(1+|\log\tfrac{|\xi|}{|\eta|}|))\big)^{-\theta}$ for
$\theta>1$.

We find more convenient to work with the martingale
difference operators $\Delta_k $ associated with the $\sigma-$algebra of all
dyadic cubes of size $2^k$ in $\mathbb R^n$ and later transfer our
results to the Littlewood-Paley operators  $\widetilde\Delta_k $.
This point of view is introduced in the next section.

We end this article with a short discussion on paraproducts, see
section \ref{s-paraproducts}. These are operators of the type
(\ref{ooo}) for specific sequences $(a_{jk})_{j,k}$ of zeros and ones.

\section{A maximal operator}\label{s-maximal}
\setcounter{equation}0

Let $(\Omega,\Sigma,\mathbb P)$ be any probability space and let
$(\Sigma_k)_{k\ge 0}$ be a {\it filtration} i.e. an increasing sequence
of sub-$\sigma-$algebras
of $\Sigma$.  We say that $(\Sigma_k)$ is a {\it dyadic filtration} if
each $\Sigma_k$ is atomic and has precisely $2^k$ atoms each with
probability $2^{-k}.$  We say $(\Sigma_k)$ is a {\it $2^n-$adic
filtration} if each $\Sigma_k$ is atomic with precisely $2^{nk}$ atoms
each with probability $2^{-nk}.$

Associated to $\Sigma_k$ we define the conditional expectation operators
$\mathcal E_kf=\mathbb E(f|\Sigma_k)$ and the martingale difference
operators $\Delta_kf =\mathcal E_kf-\mathcal E_{k-1}f$ for $k\ge 1,$ and
$f\in L_1(\Omega).$

Let $A=(a_{jk})$ be a complex  $M\times N$ matrix, and let
$(\Omega,\Sigma,\mathbb P)$ be a probability space with a dyadic
filtration $(\Sigma_k)_{k\ge 0}$. For
$1\le p<\infty$
we define $h_{p}(A)$ to be the least constant so that for all
$f\in L_p(\Omega)$  we have
\begin{equation}\label{defh}
\big\|\max_{1\le j\le M}|\sum_{k=1}^Na_{jk}\Delta_kf| \big\|_{L_p}\le
h_{p}(A)\|f\|_{L_p}.\end{equation}
 We also define the corresponding
weak constants, i.e. the least constants so that
for all $f\in L_p(\Omega)$  we have
\begin{equation}\label{defhweak}
\big\|\max_{1\le j\le M}|\sum_{k=1}^Na_{jk}\Delta_kf|\big\|_{L_{p,\infty}}\le
h_p^w(A)\|f\|_{L_p}.
\end{equation}
Finally for $0<q<p<\infty$ we define the mixed constants $h_{p,q}(A)$  as  the
least constants such that for all $f\in L_p(\Omega)$  we have
\begin{equation}\label{defhmixed}
\big\|\max_{1\le j\le M}|\sum_{k=1}^Na_{jk}\Delta_kf|\big\|_{L_q}\le
h_{p,q}(A)\|f\|_{L_p}.
\end{equation}

Note that  these definitions are independent of the choice of the probability
space and of the dyadic filtration. Indeed if $A$ is fixed, it suffices to
take $f\in L_p(\Sigma_N)$ and hence we can consider a finite
probability space with $2^N$ points and a finite dyadic filtration
$(\Sigma_k)_{k=0}^N.$
We also note that $h_p(A)$ is the operator norm of  the map
$T_A:L_p(\Omega)\to L_p(\Omega;\ell_{\infty}^M)$ defined by
$$
T_Af= \big(\sum_{k=1}^Na_{jk}\Delta_k f \, \big)_{j=1}^M.
$$
Similarly $h_p^w(A)$
is the norm of the operator $T_A:L_p\to
L_{p,\infty}(\Omega;\ell_{\infty}^M).$

Our first result is that all these constants are mutually equivalent,
when $1<p<\infty$:

\begin{Thm}\label{equiv}  If $ 1<p,q<\infty$ then there is a constant
$0<C=C(p,q)<\infty$ such that for all complex $M\times N$
matrices $A$ we have
$$
\frac1C h_p(A)\le h_q^w(A)\le h_q(A)\le Ch_p(A).
$$
\end{Thm}

\begin{proof}  It suffices to prove an estimate of the type
$h_p(A)\le Ch_q^w(A)$ for any choice of $1<p,q<\infty.$  We first prove a
weak type $(1,1)$ estimate for $T_A$, i.e. that $h_1^w(A)\le Ch_q^w(A).$
Suppose $f\in L_1$ with
$\|f\|_{L_1}=1.$ Then if $\lambda, \gamma >0,$ with $\lambda\gamma>1,$ we can
use an appropriate  Calder\'on-Zygmund decomposition 
to find   finite sets $D_1,\ldots ,D_m$ so that each
$D_l$ is an atom of some $\Sigma_l , $
$$
\gamma \lambda  \le  
 \mathbb P (D_l) ^{-1}
\int_{D_l} |f| \, d \mathbb P =\Ave_{D_l}f \le
2 \gamma  \lambda,
$$
and $|f(\omega)| \le \gamma\lambda$ if $\omega\notin \cup_{l=1}^mD_l.$
Let
$$
g=\sum_{l=1}^m (\Ave_{D_l}f)\chi_{D_l}
$$ and $E=\cup_{l=1}^mD_l.$
Then $T_A(f\chi_E-g)$ is supported in $E$ and thus
\begin{equation}\label{lou0}
\mathbb
P(\|T_A(f\chi_E-g)\|_{\ell_\infty^M}>\lambda/2) \le \mathbb P(E)\le
(\gamma\lambda)^{-1}.
\end{equation}
On the other hand $\|f-f\chi_E+g\|_{L_{\infty}}\le 3\gamma\lambda$ and
$\|f-f\chi_E+g\|_{L_1} \le 1$. Hence
$\|f-f\chi_E+g\|_{L_q} \le 3^{1/q'}(\gamma\lambda)^{1/q'}$ and so
\begin{equation}\label{lou1}
\|T_A(f-f\chi_E+g)\|_{L_{q,\infty}(\ell_\infty^M)}\le h_q^w(A)3^{1/q'}
(\gamma\lambda) ^{1/q'},
\end{equation}
which implies that
\begin{equation}\label{lou2}
\mathbb P(\|T_A(f-f\chi_E+g)\|_{\ell_\infty^M}> \lambda /2) \le
 \frac{2^q}{ \lambda^q} ( h_q^w(A) )^q  3^{q-1} (\gamma \lambda)^{q-1}  .
\end{equation}
Selecting $\gamma=1/h^w_q(A)$ and  combining with (\ref{lou0}) we obtain
(for $\lambda>h_q^w(A))$
\begin{equation}\label{lou3}
\lambda \, \mathbb P(\|T_Af\|_{\ell_\infty^M}>\lambda)\le Ch_{q}^w(A)
\end{equation}
where $C=C(p,q)$. This gives 
the weak-type (1,1) estimate for $T_A.$  Now by the Marcinkiewicz
interpolation theorem (applied to the sublinear map
$f\mapsto \|T_Af(\omega)\|_{\ell_\infty^M} $) we obtain that $h_p(A)\le
C(p,q)h_q^w(A)$ as long as $1<p<q.$

We now prove that $h_p(A)\le C(p,q)h_q^w(A)$ when $1<q<p<\infty$.
We consider the dual map
$T_A^*:L_1(\Omega;\ell_1^M)\to L_1$ defined by
$$
 T_A^*\mathbf f= \sum_{j=1}^M \sum_{k=1}^N a_{jk}\Delta_kf_j
$$ 
where $\mathbf
f(\omega)=(f_j(\omega))_{j=1}^M.$  We   have that
$T_A^*:L_r(\Omega;\ell_1^M)\to L_r$ has norm bounded by $C(q,r)h_q^w(A)$
as long as $1<r'<q$ i.e. $q'<r<\infty.$ Using this $r$ as a starting
point, we
repeat  the  argument above to show that $T_A^*:L_1(\Omega;\ell_1^M)\to
L_{1,\infty}$ has norm bounded by $Ch_q^w(A).$  The Marcinkiewicz
interpolation theorem can again be used to show that
$T_A^*:L_{p'}(\Omega,\ell_1^M)\to L_{p'}$ has norm bounded by $ Ch_q^w(A)$
for all $1<p'<r$, and thus in particular when $1<p'<q'$. Therefore we
obtain that
 $h_p(A)\le Ch_q^w(A)$ when $1<q<p<\infty$. 
\end{proof}

\noindent {\bf Remark.}  From now we will write $h(A)=h_2(A)$ so that each
$h_p(A)$ for $1<p<\infty$ is equivalent to $h(A).$

 It is of some interest to observe that even the corresponding mixed
constants are also equivalent to $h(A).$

 \begin{Thm}\label{mixed} Suppose $0<q<p$ and $1<p<\infty$.  Then there
is a constant $C=C(p,q)$ so that
$$ \frac1C h(A)\le h_{p,q}(A)\le  Ch(A).$$\end{Thm}

\begin{proof} This will depend on the following Lemma:

\begin{Lem}\label{Lem1.3} Suppose $1\le p<\infty$ and $0<q<p.$  Then there is
a constant $C=C(p,q)$ so that if
$r=\min(p,2)$ we have
\begin{equation}\label{wte}
 \|T_A\|_{L_p\to L_{r,\infty}(\ell_{\infty}^M)}\le C h_{p,q}(A).
\end{equation}
\end{Lem}

\begin{proof} (Lemma \ref{Lem1.3}) We may assume $q<r.$ This is a fairly
standard application of Nikishin's theorem, see \cite{N}.  Here we use a
version given in
\cite{Pisier}. It is simplest to consider the case when $\Omega$ is finite
with $|\Omega|=2^N.$ Consider the map
$T_A:L_p\to L_q(\Omega;\ell_{\infty}^M)$. For each
$f\in L_p$ with
$\|f\|_{L_p}\le 1$, let $F_f(x)=\|T_Af(x)\|_{\ell_{\infty}^M}.$
For  $\|f_j\|_{L_p}\le 1$ with $1\le j\le J$, $\sum_{j=1}^J|b_j|^r=1$, and
$(\epsilon_j)_{j=1}^{J}$  a sequence of independent Bernoulli random
variables on some probability space, we have
\begin{equation*}
\big\|\max_{1\le j\le J} |b_j|F_{f_j}\big\|_{L_q}  \le \mathbb
E \bigg(\big\|\sum_{j=1}^J\epsilon_jb_jT_Af_j
\big\|_{L_q(\ell_{\infty}^M)} \bigg)
 \le Ch_{p,q}(A),
\end{equation*}
since $L_p$ has type $r$.
It follows from \cite{Pisier} that there is a function $w\in L_1$, with
$\int w\,d\mathbb P=1,$ and
 $w\ge 0$ a.e such that for any set $E\subset\Omega$
$$
\Big(\int_E F_f^qd\mathbb P\Big)^{\frac1q} \le Ch_{p,q}(A) \Big(\int_E
w\, d\mathbb P\Big)^{\frac1q-\frac1r}.
$$
Now consider the set $S$ of all permutations of $\Omega$
which induce permutations of the atoms of each $\Sigma_k$ for $1\le k\le
N$; there are $2^{2^N-1}$ such permutations $\varphi.$  For $\varphi\in
S$ we have
$$
 \Big(\int_E F_{f\circ\varphi}^q \, d\mathbb P\Big)^{\frac1q} \le
Ch_{p,q}(A)\Big(\int_E w\, d\mathbb P\Big)^{\frac1q-\frac1r}
$$
or equivalently
$$
 \Big(\int_E F_{f}^q \, d\mathbb P\Big)^{\frac1q} \le Ch_{p,q}(A)\Big(\int_E
w\circ\varphi^{-1}\, d\mathbb P\Big)^{\frac1q-\frac1r}.
$$
Raising to the power $(\frac1q-\frac1r)^{-1}$, averaging over $S$, and then
raising to the power $\frac1q-\frac1r$ gives
$$
\Big(\int_E F_{f}^q d\mathbb P\Big)^{\frac1q} \le
Ch_{p,q}(A)\bigg(\frac{1}{|S|}\sum_{\varphi\in S}\Big(\int_E
w\circ\varphi^{-1} \,  d\mathbb P \Big) \bigg)^{\frac1q-\frac1r}.
$$
But this implies
$$
 \Big(\int_E F_{f}^q d\mathbb P\Big)^{\frac1q} \le
Ch_{p,q}(A)\, \mathbb P(E)^{\frac1q-\frac1r}
$$
which gives the required weak type estimate (\ref{wte}).
\end{proof}

We now return to the proof of Theorem \ref{mixed}.
We first observe that we always have $h_{p,q}(A)\le Ch_p^w(A)$ since
$q<p$.
If $1<p\le 2$, Lemma \ref{Lem1.3} gives that
$h_p^w(A)\le Ch_{p,q}(A)$ and the required conclusion follows from Theorem
\ref{equiv}.  Assume therefore that $p>2$ and that $T_A$ maps
$L_p\to L_q(\ell_{\infty}^M) $ with norm $h_{p,q}(A)$.
Fix $f$ with $\|f\|_{L_1}=1$ and use the
Calder\'on-Zygmund decomposition  of Theorem \ref{equiv},
to obtain   (\ref{lou0})  as before, but instead of (\ref{lou1}) the
estimate
\begin{equation}\label{lou1-1}
\|T_A(f-f\chi_E+g)\|_{L_q }\le h_{p,q} (A) 3^{1/p'}
(\gamma\lambda) ^{1/p'},
\end{equation}
which implies
\begin{equation}\label{lou2-2}
\mathbb P(\|T_A(f-f\chi_E+g)\|_{\ell_\infty^M}> \lambda /2) \le
 \frac{2^q}{ \lambda^q} ( h_{p,q} (A) )^q  3^{q/p'} (\gamma \lambda)^{q/p'}.
\end{equation}
Selecting $\gamma= h_{p,q}(A)^{-s}\lambda^{s-1}$ with
$\frac1s=\frac1{p'}+\frac1q$ and  combining with
(\ref{lou0}) we obtain
\begin{equation}\label{lou3-3}
\lambda \, \mathbb P(\|T_Af\|_{\ell_\infty^M}>
\lambda)^{\frac1s}\le C h_{p,q}(A). 
\end{equation}
This says that $T_A$ maps $L_1$ into $ L_{s,\infty}(\ell_{\infty}^M)$   with
norm at most
$Ch_{p,q}(A)$,  in particular that $T_A$ maps $L_1$ into  $L_t(\ell_{\infty}^M)$
as long as $0<t<s$.     Lemma \ref{Lem1.3}   gives that
$T_A$ maps $L_p $ into $L_{2,\infty}(\ell_{\infty}^M)$ and also
$L_1 $ into $L_{2,\infty}(\ell_{\infty}^M)$ with norms at most a
multiple of $h_{p,q}(A)$. By interpolation it follows that
$T_A$ maps $L_r $ into $L_{2,\infty}(\ell_{\infty}^M) \subset
L_{r,\infty}(\ell_{\infty}^M)$ for $1\le r\le 2.$
We conclude that
$h_r^w(A)\le Ch_{p,q}(A)$ for $1< r< 2$ but since
$h_r^w(A)$ is comparable to $h_p^w(A)$, we finally obtain
$h_p^w(A)\le Ch_{p,q}(A)$.  Since the converse inequality is always
valid when $q<p$, we apply Theorem \ref{equiv} to conclude  the proof.
\end{proof}

We next prove the elementary observation for $1<p<\infty,$ that $h(A)$
remains unchanged when interpolating extra columns   or
extra rows of zeros.

\begin{Lem}\label{zeros}  Let $A$ be a complex $M\times N$ matrix  and
$(m_r)_{r=1}^M,
(n_s)_{s=1}^N$ be two
increasing finite sequences of natural numbers.
Suppose $M_1\ge m_M$ and $N_1\ge n_N.$
Let $B=(b_{jk})$ be the $M_1\times N_1$-matrix
defined by $b_{jk}=a_{rs}$ when $j=m_r$ and $k=n_s$, and $b_{jk}=0$
otherwise.  Then $h(A)= h(B).$
\end{Lem}

\begin{proof} Interpolating extra rows of zeros is trivial, so we can
assume $m_r=m$ for all $r.$ For the case of columns, we    only
need to show that  $h(B)\le h(A).$
We may suppose that
$\Omega$ is a finite set with $2^{N_1}$ points and that
$(\Sigma_k)_{k=0}^{N_1}$ is a finite dyadic filtration of $\Omega$.  It is then
possible to write $\Omega=\Omega_1\times \Omega_2$ where
$|\Omega_1|=2^{N_1-N }$ and $|\Omega_2|=2^{N  }$,
  and find a dyadic filtration
$(\Sigma_k^{(1)})_{k=0}^{N_1-N}$ of $\Omega_1$ and a dyadic  filtration
$(\Sigma_k^{(2)})_{k=0}^{N}$ of $\Omega_2$ so that
$\Sigma_k^{(1)}\times\Sigma_k^{(2)}=\Sigma_{n_k},$ for $0\le k\le N$ and
$\Sigma_{k+1}^{(1)}\times \Sigma_k^{(2)}=\Sigma_{n_{k+1}-1}$ for $0\le
k\le N-1.$ Then for
$f\in L_2(\Omega_1\times\Omega_2)$ let $g=\sum_{k=1}^N \Delta_kf$ and
note that
$$
\Delta_{n_k}f(\omega_1,\omega_2)=\Delta_k^{(2)}g_{\omega_1}(\omega_2), 
$$
where $g_{\omega_1}(\omega_2)=g(\omega_1,\omega_2).$  Hence
$$
\int_{\Omega_2}
\sup_j \Big|\sum_{k=1}^Na_{j,n_k}\Delta_{n_k}f(\omega_1,\omega_2)
\Big|^2d\omega_2
\le h_p(A)\int_{\Omega_2}|g(\omega_1,\omega_2)|^2d\omega_2.
$$
Integrating over $\Omega_1$ gives
$$
\big\|\sup_j|\sum_{k=1}^Na_{j,n_k}\Delta_{n_k}f|
\big\|_{L_2} \le h_p(A)\|g\|_{L_2}\le
h_p(A)\|f\|_{L_2}.
$$
This completes the proof.\end{proof}

We can now    extend our definitions, replacing   dyadic filtrations by
$2^n$-adic filtrations:

\begin{Prop}\label{filtration}
Suppose $n\in\mathbb N$ and $1<p<\infty$.
Then there is a constant $C(p,n)$ with the following property. Let
$(\Omega,\Sigma,\mathbb P)$ be a probability space and suppose
$(\Sigma_k)_{k=0}^{\infty}$ is a $2^n$-adic filtration. Let
$A$ be any
$M\times N$ matrix and let $h_p(A;n)$ be the least constant so that
$$
\big\|\sup_j|\sum_{k=1}^Na_{jk}\Delta_kf|\big\|_{L_p}\le h_p(A;n)\|f\|_{L_p} ,
$$
and $h_p^w(A;n)$ be the least constant so that
$$
\big\|\sup_j|\sum_{k=1}^Na_{jk}\Delta_kf|\big\|_{L_{p,\infty}}\le
h^w_p(A;n)\|f\|_{L_p}.
$$
Then $h_p^w(A)\!\le\! h_p^w(A;n),\ h_p(A)\!\le\! h_p(A;n)$,
and $h_p^w(A;n)\!\le\! h_p(A;n)\!\le\! Ch(A).$
\end{Prop}

\begin{proof} This is essentially trivial; we need only to prove that
$h_p(A;n)\le Ch(A).$  To do this note that $h_p(A;n)=h_p(B)$ where $B$ is
obtained from $A$ by repeating each column $n$ times. The proposition
follows then by the triangle law from Lemma \ref{zeros}.\end{proof}

\section{Estimates for $h(A)$}\label{estimates}
\setcounter{equation}{0}

We next turn to the problem of estimating $h(A).$  We shall assume that
$(\Omega,\mathbb P)$ is a fixed probability space with a dyadic
filtration
$(\Sigma_k)_{k=0}^{\infty}.$
 Our first estimate is trivial.

\begin{Prop}\label{bv}   There is a constant $C$
so that for any $M\times N$ matrix
$A=(a_{jk})$ we have
$$ h(A)\le C \sup_{1\le j\le M}\sum_{k=0}^N
|a_{jk}-a_{j,k+1}|,$$
where we set $a_{j0}=a_{j,N+1}=0$ for all $1\le j\le M.$
\end{Prop}

\begin{proof}  Suppose $f\in L_2.$  Summation by parts gives
$$
 \sum_{k=1}^Na_{jk}\Delta_kf = \sum_{k=0}^N
(a_{jk}-a_{j,k+1}){\mathcal E_k}f,
$$
thus
$$
|\sum_{k=1}^Na_{jk}\Delta_kf|\le (\sup_{1\le j\le M} \sum_{k=0}^N
|a_{jk}-a_{j,k+1}|)\sup_k|{\mathcal E}_kf|,
$$
and the result follows because of the maximal estimate  
$$ \|\sup_k|{\mathcal E}_kf|\|_{L_2}\le C\|f\|_{L_2},$$
proved in   \cite{garsia}. 
\end{proof}

We next turn to the problem of getting a more delicate estimate.  To do
this we need the theory of a certain Lorentz space.  Let
$w=(w_k)_{k=1}^{\infty}$ be a decreasing sequence of positive real numbers.
We will consider the following two conditions on $w:$
\begin{equation}\label{cdn1}
\exists C>0,\,\,\,\exists\theta>0,\qquad w_k\le
C\Big(\frac{\log(j+1)}{\log(k+1)}\Big)^{\theta}w_j \quad\text{when }\,\, 1\le
j\le k,\end{equation}
(where throughout this paper $\log$ denotes  the logarithm with base $2$) and
\begin{equation}\label{cdn2}
\sum_{k=1}^{\infty}\frac{w_k}{k}<\infty.\end{equation}
Roughly speaking (\ref{cdn1}) means that $w_k$ decays logarithmically
while (\ref{cdn2}) implies that it decays reasonably fast.  Note that
$w_k=(\log(k+1))^{-\theta}$ satisfies (\ref{cdn1}) if $\theta>0$ and
(\ref{cdn2}) if $\theta>1.$  The sequence
$w_k=(\log(k+1))^{-1}(\log\log(k+2))^{-\theta}$ satisfies
both (\ref{cdn1}) and  (\ref{cdn2}) when $\theta>1.$

Now let $d=d(w,1)$ be the Lorentz sequence space of all complex sequences
$\mathbf u=(u_k)_{k\in\mathbb Z}$ such that
$$\|\mathbf u\|_{d}=\sup_{\pi}\sum_{k\in\mathbb Z}w_{\pi(k)}|u_k|
<\infty$$ where the supremum is taken over all one-one maps $\pi:\mathbb
Z \to\mathbb N.$  The dual of $d(w,1)$ can be naturally identified as the
space $d^*=d^*(w,1)$ consisting of all sequences $(v_k)_{k\in\mathbb Z}$
so that
$$
 \sup_{k\in\mathbb N}
\frac{ v_1^*+\cdots+v_k^* }{w_1+\cdots+w_k}=\|\mathbf
v\|_{d^*}<\infty
$$
where $(v_k^*)_{k=1}^{\infty}$ is the decreasing rearrangement of
$(|v_k|)_{k\in\mathbb Z}.$  We refer to \cite{LT1} p. 175 for properties of
Lorentz spaces. Note that under condition (\ref{cdn1}), $d(w,1)$ is also an
Orlicz sequence space (see \cite{LT1} p. 176).

The following Lemma is surely well-known to specialists, but we include a
proof.

\begin{Lem}\label{cotype}
Under condition (\ref{cdn1}), the Lorentz space $d(w,1)$ has cotype two.
\end{Lem}

\begin{proof} By combining Proposition 1.f.3 (p.82) and Theorem 1.f.7
(p.84) of \cite{LT2} one sees that it is only necessary to show that
$d(w,1)$ has a lower $q$-estimate for some $q<2.$  To do this observe
that if $\mathbf v_1,\cdots,\mathbf v_N$ are disjointly supported
sequences, then
$$ \|\sum_{j=1}^N\mathbf v_j\|_{d}\ge
\inf_{k\ge 1}\frac{w_k}{w_{kN}}\sum_{j=1}^N\|\mathbf v_j\|_{d}.$$
Hence
$$ \sum_{j=1}^N\|\mathbf v_j\|_{d(w,1)}\le C(\log
(N+1))^{\theta}\|\sum_{j=1}^N\mathbf v_j\|_{d}.$$
Now suppose $1<q<2$ and $\|\sum\mathbf v_j\|_d=1.$
Then for each $s\in \mathbb N$, let $m_s$ be the number of $j$ so that
$2^{-s}<\|\mathbf v_k\|_d \le  2^{-s+1} .$
Then
$$ m_s2^{-s} \le C(\log (m_s+1))^{\theta}.$$
This in turn implies that
$$ m_s^{1-\rho} \le C2^s$$ where $\rho>0$ is chosen so that
$(1-\rho)^{-1}<q.$  Then we obtain an estimate
$$ \sum_{j=1}^N\|\mathbf v_j\|_d^q\le C\sum_{s=1}^{\infty}m_s2^{-sq}\le
C'.$$  This establishes a lower $q$-estimate.\end{proof}

The norms $\|\cdot\|_d$ and $\|\cdot\|_{d^*}$ are  of course
defined for finite sequences with $M$ elements and  thus can be thought as
norms on $\Bbb C^M.$ We denote these spaces
$d(w,1)^{(M)}$ and $d^*(w,1)^{(M)}.$

\begin{Prop}\label{changesign}
If $(w_n)$ satisfies both
(\ref{cdn1}) and (\ref{cdn2}) then given $2<p<\infty$ there is a
constant
$C$ so that for any sequence $\epsilon_k=\pm1$ and any $M,N\in\mathbb N$
we have the estimate
$$
\Big(\mathbb E\big( \big\|\sum_{k=1}^N\epsilon_k\Delta_k\mathbf
f\big\|_{\ell_{\infty}}^2\big) \Big)^{\frac12}\le C
 \big(\mathbb E \big(\|\mathbf f\|_{d^*} ^p\big)
\big)^{\frac1p},
$$
for any $\mathbf f\in L_p(\Omega;d^*(w,1)^{(M)}).$
\end{Prop}

\begin{proof}  We start by using an argument due to Muckenhoupt
\cite{muckenhoupt}, see also \cite{Torchinsky}.  For any fixed
$\epsilon_1,\ldots,\epsilon_N$ let $S=\sum_{k=1}^N\epsilon_k\Delta_k.$
Now fix $f\in L_{\infty}.$
Then by a result of Burkholder \cite{burkholder}, $\|S\|_{L_p\to L_p}=p-1$ if
$2\le p<\infty.$    Then for any $\alpha>0$ we
have
\begin{equation}\label{ghg}
\mathbb E(\cosh (\alpha |S f|)) \le
1+\sum_{m=1}^{\infty}\frac{\alpha^{2m}}{(2m)!}(2m-1)^{2m}\|f\|_{L_{2m}}^{2m}.
\end{equation}
Since $\|f\|_{L_{2m}}^{2m}\le \|f\|_{L_2}^2 
\|f\|_{L_{\infty}}^{2m-2}$ and since for
$m\ge 1$ we have
$$
 \frac{(2m-1)^{2m}}{(2m)!} \le  \frac{(2m )^{2m}}{(2m)!} \le e^{2m} ,
$$
it follows from (\ref{ghg}) that
$$ 
\mathbb E(\cosh (\alpha |S f|)-1) \le
(\alpha e)^2\|f\|_{L_2}^2\sum_{k=0}^{\infty}(\alpha e)^{2k}
\|f\|_{L_\infty}^{2k}.
$$
In particular if $\alpha e\|f\|_{\infty}\le \frac12$ we have
\begin{equation}\label{firstestimate}
 \mathbb E(\cosh (\alpha |S f|)-1) \le 2e^2
\alpha^2\|f\|_{L_2}^2.
\end{equation}

At this point we return to the Lorentz space $d(w,1)$. Let us define
$\gamma_0=0$,
$\gamma_1=1$,  and   $\gamma_k=2^{2^{k-2}}$ for $k\ge 2.$  Let
$W_k=w_{\gamma_k}.$
It will be convenient to normalize condition (\ref{cdn2}) so that we have
\begin{equation}\label{cdn20}
\sum_{k=1}^{\infty}\gamma_{k}W_k=1.
\end{equation}
We also note that
(\ref{cdn1})
implies the existence of a constant $C$ so that we have
\begin{equation}\label{cdn10}|w_1+\cdots+w_k|\le Ckw_k\end{equation} for
$k\ge 1.$

Now suppose $\mathbf f=(f_j)_{j=1}^M\in L_{\infty}(\Omega;\mathbb C^M).$
Suppose that $\mathbf f$ is supported on a measurable set $E$ and satisfies
$\|\mathbf f(\omega)\|_{d^*}\le 1$ everywhere.  Then we can define
a measurable map $\pi$ from $\Omega$ into the set of permutations of
$\{1,2,\ldots,M\}$ so that $|f_{\pi(\omega)(1)}(\omega)|\ge
|f_{\pi(\omega)(2)}(\omega)|\ge\cdots \ge |f_{\pi(\omega)(M)}(\omega)|$ for all
$\omega \in  \Omega .$ Thus
$$ |f_{\pi(\omega)(j)}(\omega)|\le Cw_j$$ for all $1\le j\le M.$  Let
$E_{jk}=\{\omega\in E:\ \pi(\omega)(k)=j\} $ when $j,k\in \{1,\dots , M\}$ and
$E_{jk}=\emptyset$ otherwise.
Now for $1\le j\le M$ and $l=1,2,3,\dots$, let
$$
f_j^{(l)}=\sum_{k=\gamma_{l-1}}^{\gamma_l-1}f_j\chi_{E_{jk}}
$$
so that $f_j=\sum_{l=1}^\infty f_j^{(l)}$.
If $0<\alpha e\le \frac1{2C }$ we can estimate
\begin{equation*}\begin{align*}
\mathbb E(\cosh (\alpha |S f_j|)-1)
 &=\mathbb E
\big(\cosh \big(\big|\sum_{l=1}^{\infty}\alpha Sf_j^{(l)}\big|\big)-1\big)\\
&\le \mathbb E\Big(\max_{l\ge 1}\big(\cosh (\alpha
\gamma_l^{-1}W_l^{-1}|Sf_j^{(l)}| )-1\big)\Big)\\ &\le  e^2\alpha ^2
\sum_{l=1}^{\infty}\gamma_l^{-2}W_l^{-2}\|f_j^{(l)}\|_{L_2}^2,
\end{align*}\end{equation*}
in view of (\ref{firstestimate}) since
 $\|f_j^{(l)}\|_{L_\infty}\le CW_l$ and  $\alpha
\gamma_l^{-1}W_l^{-1}\|f_j^{(l)}\|_{L_\infty}\le \frac12. $ Thus
$$
\mathbb E(\cosh (\alpha |S f_j|)-1) \le
e^2 C^2\alpha^2\sum_{l=1}^{\infty}\gamma_l^{-2}
\sum_{k=\gamma_{l-1}}^{\gamma_l-1} \mathbb P(E_{jk}).
$$
It follows that
$$
\mathbb E(\cosh (\alpha \|S\mathbf f\|_{\ell_{\infty}})-1) \le
e^2C^2\alpha^2\sum_{j=1}^M\sum_{l=1}^{\infty}\gamma_l^{-2}
\sum_{k=\gamma_{l-1}}^{\gamma_l-1} \mathbb P(E_{jk}).
$$
Note that for each $k\in\mathbb N,$ $\sum_{j=1}^M\mathbb
P(E_{jk})\le
\mathbb P(E).$ Hence we obtain that if $\mathbf f$ is supported on
$E$ with $\|\mathbf f(\omega)\|_{d^*}\le 1$ everywhere  and
$\alpha e<\frac1{2C},$ then
\begin{equation}\label{firststep}
 \mathbb E(\cosh (\alpha \|S\mathbf
f\|_{\ell_{\infty}})-1) \le e^2 C^2\alpha^2
\sum_{l=1}^{\infty}\gamma_l^{-1}\mathbb P(E)=C_1\alpha^2\mathbb
P(E)\end{equation}  for a suitable constant
$C_1.$
Let us next refine (\ref{firststep}).  For $n\ge 0,$ let
$$
G_n=\{\omega\in E: \ 4^{-n-1}< \|\mathbf(\omega)\|_{d^*}\le 4^{-n}\}.
$$
Then by (\ref{firststep}) we have if $\alpha<(4Ce)^{-1}$
$$
 \mathbb E(\cosh (2^{n+1}\alpha \|S (\mathbf
f\chi_{G_n})\|_{\ell_{\infty}})-1) \le C_1\alpha^2 4^{-n} \mathbb P(G_n)
$$
and as
$$
\mathbb E(\cosh (\alpha \|S\mathbf f\|_{\ell_{\infty}})-1)\le
\mathbb E\Big(\sup_{n\ge 0}\big(\cosh (2^{n+1}\alpha\| S(\mathbf
f\chi_{G_n})\|_{\ell_{\infty}})-1\big)\Big),
$$
we obtain, under the assumptions
$\|\mathbf f(\omega)\|_{d^*}\!\le\! 1$ everywhere and $\alpha\!<\! (4C)^{-1}$,
\begin{equation}\label{2step}
 \mathbb E(\cosh (\alpha \|S\mathbf f\|_{\ell_{\infty}})-1)\le
C_1\alpha^2\sum_{n=0}^{\infty}4^{-n}\mathbb P(G_n)\le C_2\mathbb
E(\|\mathbf f\|_{d^*}).\end{equation}

If we use a fixed value of $\alpha$ and the estimate $x^2\le 2(\cosh
x-1)$ we find that
$$
 \mathbb E(\|S\mathbf f\|_{\ell_{\infty}}^2)\le C_3\mathbb E (\|\mathbf
f\|_{d^*})
$$
if $\|\|\mathbf f\|_{d^*}\|_{\infty}\le 1$.  This in turn gives
us for every $\mathbf f\in L_{\infty}(\Omega; d^*(w,1)^{(M)})$
\begin{equation}\label{3step}
\mathbb E(\|S\mathbf f\|_{\ell_{\infty}}^2)\le C_3 \|\|\mathbf
f\|_{d^*}\|_{\infty}\, \mathbb E(\|\mathbf f\|_{d^*}).\end{equation}

Now let $2<p<\infty$ and fix $\mathbf f$ with
$\mathbb E(\|\mathbf f\|_{d^*}^p)= 1$. We set
$E_0=\{\|\mathbf f\|_{d^*}\le 1\}$ and $E_n=\{2^{n-1}<\|\mathbf
f\|_{d^*}\le 2^n\}$ for $n\ge 1.$  Applying (\ref{3step}) we obtain
 \begin{align*}
(\mathbb E(\|S\mathbf f\|_{\ell_{\infty}}^2))^{\frac12}  \le  &
\big( C_3 \sum_{n=0}^{\infty} 2^n
\mathbb P(E_n) \mathbb E ( \|\mathbf f\|_{d^*} )
\big)^{\frac12}\le
 C_3^{\frac12}\sum_{n=0}^{\infty}2^{\frac{n}{2}}\mathbb P(E_n)^{\frac12} \\
\le  & C_3^{\frac12}\big(\sum_{n=0}^{\infty}
2^{(2-p)n} \big)^{\frac12}
\big(\sum_{n=0}^{\infty} 2^{np}\mathbb
P(E_n)  \big)^{\frac12} \le C_4 ,
 \end{align*}
which completes the proof under the assumption
$\mathbb E(\|\mathbf f\|_{d^*}^p)= 1$. The general case follows by scaling.
\end{proof}

We now establish our main estimate for $h(A).$

\begin{Thm}\label{lorentz}
Let $w=(w_n)_{n=1}^{\infty}$ be a sequence
satisfying
(\ref{cdn1}) and (\ref{cdn2}).  Then there is a constant $C$ so that
for any $M\times N$ matrix $A=(a_{kj})_{j,k}$ we have
 $$ h(A) \le C \max_{1\le k\le N}\|\mathbf a_k\|_{d^*}$$ where
$\mathbf a_k=(a_{kj})_{j=1}^M.$  In particular we have
$$ h(A) \le C\max_{j,k}\frac{|a_{jk}|}{w_{|j-k|+1}}.$$\end{Thm}

\begin{proof}  We suppose $p>2$ and that $A$ is a matrix satisfying
$\max_{1\le k\le N}\|\mathbf a_k\|_{d^*}\le 1.$  Consider the
operator
$T_A:L_p(\Omega)\to L_2(\Omega;\ell_{\infty}^M).$  The adjoint operator
is $T_A^*:\,\, L_2(\Omega;\ell_1^M)\to L_{p'}(\Omega)$  given by
$$
T_A^*(\mathbf f)=\sum_{k=1}^N\langle \Delta_k\mathbf f,\mathbf
a_k\rangle.
$$
The dual statement of the result in Proposition
\ref{changesign} gives that for any sequence of $\pm 1$'s,
$\epsilon_1,\ldots,\epsilon_N$ we have the estimate
\begin{equation}\label{signs}
\big(\mathbb E \big(\|\sum_{k=1}^N\epsilon_k\Delta_k\mathbf
f\|_{d}^{p'}\big)\big)^{\frac1{p'}}\le
C(\mathbb E(\|\mathbf f\|_{\ell_1}^2))^{\frac12}  \end{equation}
where $C$ depends only on $(w_n).$  Now let
$\epsilon_1,\ldots\epsilon_N$ be a sequence of independent Bernoulli
random variables on some probability space $(\Omega',\mathbb P').$  We
use $\mathbb E'$ to denote expectations on $\Omega'.$
Using Lemma \ref{cotype} we obtain
 \begin{align*}
(\mathbb E (\|T_A^*\mathbf f\|_d^{p'}))^{\frac1{p'}} 
&\le C_0\big(\mathbb
E \big(\sum_{k=1}^N |\langle \Delta_k\mathbf f,\mathbf
a_k\rangle|^2 \big)^{\frac {p'}2}\big)^{\frac1{p'}}\\
&\le C_0 \big(\mathbb E \big(\sum_{k=1}^N\|\Delta_k\mathbf
f\|_d^2\big)^{\frac {p'}2}\big)^{\frac1{p'}} \\
&\le C_1 \big(\mathbb E\mathbb E'\big( \|\sum_{k=1}^N\epsilon_k\Delta_k\mathbf
f\|_d^{p'} \big)\big)^{\frac1{p'}}\\
&\le C_2 (\mathbb E\|\mathbf
f\|_{\ell_1}^{2})^{\frac1{2}}.
\end{align*}
This gives $h_{p,2}(A)\le C_2$ which completes the proof by using
Theorem \ref{mixed}.
\end{proof}

\noindent {\bf Remark.}  Theorem \ref{lorentz} implies that given any
$\theta>1$ there is a constant $C_{\theta}$ so that
\begin{equation}\label{44444}
h(A)\le C_{\theta}
\end{equation}
whenever $A =(a_{kj})_{j,k}$ is a
matrix satisfying
\begin{equation}\label{logcondition}
|a_{jk}|\le 2 (\log (2+|j-k|))^{-\theta}.
\end{equation}
We show that this is not the case when $0<\theta<\frac12$.
Let $N$ be any natural number and define $A=(a_{jk})$ to be
a $2^N\times N$ matrix given by  $a_{jk}=b_{jk}N^{-\theta}$, where
$b_{jk}= \pm1 $
and the set
$(b_{jk})_{j=1}^{2^N}$ runs through all $2^N$ choices of signs.  Choose $f$
real so that
$|\Delta_k
f|=1$ for $1\le k\le N.$  Then $\|f\|_{L_2}=\sqrt N.$  On the other hand
$$
 \max_{1\le j\le 2^N}|\sum_{k=1}^{N}a_{jk}\Delta_kf|=N^{1-\theta}\,
\chi_{\Omega},
$$
which implies that  $h(A)\ge   N^{\frac12-\theta}$. However
$$
|a_{jk}|\le N^{-\theta} \le 2(N+1)^{-\theta} \le 2 (\log (2+|j-k| ))^{-\theta}
$$
but   $h(A)\ge   N^{\frac12-\theta}\to \infty$ as $N\to \infty$.
Thus (\ref{44444}) fails when $0<\theta<\frac12$.

\section{The harmonic version of the maximal operator}\label{s-harmonic}
 \setcounter{equation}0

We shall now fix $n\in\mathbb N$ and work with $\mathbb R^n.$
Let $\mathcal D_0$ be the collection of all unit cubes of the form
$\prod_{j=1}^n[m_j,m_j+1]$ where $m_j\in\mathbb Z$ and let $\mathcal D_k$
be the set of all cubes of the form
$\prod_{j=1}^n[2^{-k}m_j,2^{-k}(m_j+1)]$ where
$m_j\in\mathbb Z.$ For
$k\in \mathbb Z$, let $\Sigma_k$ denote  the $\sigma-$algebra generated by the
dyadic cubes
$\mathcal D_k$.  We define the corresponding conditional expectation
operators
$$
 \mathcal E_kf=\sum_{Q\in \mathcal D_k}(\Ave_{Q}f)\chi_Q
$$
for $f\in
L_1^{loc}(\mathbb R^n)$ and the martingale difference operators
$\Delta_kf=\mathcal E_kf-\mathcal E_{k-1}f $ for $k\in \mathbb Z$.

Now let $A=(a_{jk})_{j,k\in\mathbb Z}$ be any infinite complex
matrix.  We shall call $A$   a $c_{00}-$matrix if it has only finitely many
non-zero entries. For a $c_{00}-$matrix define
$h_p[A;n]$ to be the least constant such that for all $f\in L_p(\mathbb R^n)$
we have
\begin{equation}\label{defh2}
\|\max_{j\in\mathbb Z}|\sum_{k\in\mathbb Z}a_{jk}\Delta_kf|\|_{L_p}\le
h_{p}[A;n]\|f\|_{L_p}.\end{equation}
Also let $h_p^{w}[A;n]$ be the corresponding weak-type constant, i.e. the
least constant such that for all $f\in L_p(\mathbb R^n)$
we have
\begin{equation}\label{defhweak2}
\|\max_{j\in\mathbb Z}
|\sum_{k=1}^Na_{jk}\Delta_kf|\|_{L_{p,\infty}}\le
h_p^w[A;n]\|f\|_{L_p}.
\end{equation}
The following Lemma is  easily verified  and we omit its proof.
\begin{Lem}\label{lem3.1}
Let $h_p^w(A;n)$ and $h_p(A;n)$ be    as in Proposition
\ref{filtration}.
For any $1<p<\infty$  and any infinite $c_{00}$-matrix $A$ we
have
$h_p[A;n]=h_p(B;n)$ and $h_p^w[A;n]=h_p^w(B;n)$, where $B$ is any $M\times
N$ matrix of the form $b_{jk}=a_{j+r,k+s}$ for some $r,s\in\mathbb Z$
such that $a_{j+r,k+s}=0$ unless $1\le j\le M$ and $1\le k\le N$.
\end{Lem}

Now for any infinite matrix $A$ we define
$$
h(A)=\sup_Nh \big((a_{j-N,k-N})_{1\le j\le 2N}^{1\le k\le 2N}\big).
$$
The following is an immediate consequence of   Lemma \ref{lem3.1} and
Proposition \ref{filtration}.

\begin{Cor}
For any $1<p<\infty$ and any $n\in\mathbb N$  there is a
constant $C=C(p,N)$ so that for any infinite $c_{00}$-matrix we have
$$C^{-1}h(A)\le h_p^w[A;n]\le h_p[A;n]\le Ch(A).$$
\end{Cor}

We now turn to the harmonic model of the maximal operator
studied in section \ref{s-maximal}.
Let $\mathcal S(\mathbb R^n)$ denote the set of all Schwartz functions
on $\mathbb R^n$  and  for $f\in \mathcal S(\mathbb R^n)$ let
$$
\widehat{f}(\xi)= \int_{\mathbb R^n} f(x)e^{- 2\pi i\langle \xi , x\rangle} dx
$$
denote the Fourier transform of $f$.  We will denote by
$f\spcheck (\xi) =\widehat{f}(-\xi)$ the inverse Fourier transform of $f$.
 We shall fix a radial function
$\psi\in\mathcal S(\mathbb R^n)$   whose Fourier transform is real-valued and
satisfies $\widehat\psi(\xi)=1$ for $|\xi|\le 1$ and $\widehat\psi(\xi) =0$
for $|\xi|\ge 2.$
We define a Schwartz function $\phi$ by setting $\widehat{\phi}(\xi)=
\widehat{\psi}(\xi)-\widehat{\psi}(2\xi)$. Then $\widehat{\phi}$ is
supported in the annulus $2^{-1}\le |\xi |\le 2$.
 We then define
$\psi_j(x)=2^{nj}\psi(2^{j}x)$ and
$\phi_j(x)=2^{nj}\phi(2^{j}x) $ for $j\in\mathbb Z.$  Note that
$\widehat{\phi_j}(\xi)=\widehat{\phi }(2^{-j}\xi)$ is supported in the annulus
$2^{j-1}\le |\xi |\le 2^{j+1}$. We also define operators
$$
 \widetilde S_jf=\psi_j *f \quad\text{and }\,\,\,
\widetilde\Delta _j f=\phi_j*f
$$
for $f\in L_1+L_{\infty}.$  The   $\widetilde\Delta _j$'s are
called the Littlewood-Paley operators. Now if
$A=(a_{jk})_{(j,k)\in\mathbb Z^2}$ is an infinite $c_{00}$-matrix and
$1<p<\infty,$ we let $\widetilde
h_p(A)$ be the least constant so that  for all  $f\in L_p$ we have
\begin{equation}\label{tildeh}
  \big\|\sup_{j\in\mathbb Z}\big|\sum_{k\in\mathbb
Z}a_{jk}\widetilde\Delta_k f\big| \big\|_{L_p}\le \widetilde h_p(A)\|f\|_{L_p}.
\end{equation}
We also define   $\widetilde h_p^w(A)$ to be the least
constant such that for all  $f\in L_p$ we have
\begin{equation}\label{tildehw}
  \big\|\sup_{j\in\mathbb Z} \big|\sum_{k\in\mathbb
Z}a_{jk}\widetilde\Delta_k f\big| \big\|_{L_{p,\infty}}\le 
\widetilde h^w_p(A)\|f\|_{L_p}.
\end{equation}

We now have the following.

\begin{Lem}\label{translation}
Suppose $r\in\mathbb Z$.  Then if
$1<p<\infty$ and $A=(a_{jk})$ is any infinite $c_{00}$-matrix, then $\widetilde
h_p(A)=\widetilde h_p(B)$ and $\widetilde h_p^w(A)=\widetilde h_p^w(B)$, where
$B=(b_{jk})$ and $b_{jk}=a_{j,k+r}.$
\end{Lem}

\begin{proof}
Consider the dilation operator $D_rf(x)=f(2^{-r}x).$  Then
$D_r^{-1}\widetilde\Delta_kD_rf=  \widetilde\Delta_{k-r}f   $ and we have
\begin{align*}
&\big\|\sup_j\big|\sum_k a_{j,k+r}\widetilde\Delta_k f\big|  \big\|_{L_p}
 =\big\|\sup_j\big|\sum_ka_{jk}\widetilde\Delta_{k-r}f\big|  \big\|_{L_p}\\
 =&2^{-rn/p} \big\|\sup_j\big|\sum_k a_{jk}\widetilde\Delta_kD_rf\big|
\big\|_{L_p}
\le 2^{-rn/p}h_p(A)\|D_rf\|_{L_p}=h_p(A)\|f\|_{L_p},
\end{align*}
which implies $\widetilde h_p(B)\le \widetilde h_p(A).$  Likewise we obtain
$\widetilde h_p(A)\le \widetilde h_p(B).$
The  corresponding result for the weak type constants follows similarly.
\end{proof}

Next we prove  that the Littlewood-Paley
operators $\widetilde \Delta_j$ and the martingale difference
operators $\Delta_k$
are   essentially orthogonal on $L_2$ when $k\neq j$.

\begin{Prop}\label{LP-martingale}
There exists a   constant $C$ so  that
for every $k,j$ in $\mathbb Z$ we have  the following estimate
on the operator norm of $\Delta_j\widetilde \Delta_k:\,
 L_2(\mathbb R^n)\to L_2(\mathbb R^n)$
\begin{equation}\label{4444}
 \|\Delta_k\widetilde \Delta_j\|_{L_2\to L_2}\le C2^{-|j-k|}.
\end{equation}
\end{Prop}

\begin{proof}
By a simple dilation argument it suffices to prove (\ref{4444}) when $k=0$.
In this case  we have the estimate
\begin{align*}
&\|\Delta_0\widetilde \Delta_j\|_{L_2\to L_2} =
\|\mathcal E_0\widetilde \Delta_j-\mathcal E_{-1}
\widetilde \Delta_j \|_{L_2\to  L_2} \\
\le &\|\mathcal E_0\widetilde \Delta_j- \widetilde \Delta_j \|_{L_2\to
L_2} + \|\mathcal E_{-1}\widetilde \Delta_j- \widetilde \Delta_j \|_{L_2\to
L_2}
\end{align*}
and also by the self-adjointness of the $\Delta_k$'s and $\widetilde\Delta_j$'s
we have
\begin{align*}
&\|\Delta_0\widetilde \Delta_j\|_{L_2\to L_2} =\|\widetilde
\Delta_j\Delta_0 \|_{L_2\to L_2} =
\|\widetilde \Delta_j\mathcal E_0-
\widetilde \Delta_j \mathcal E_{-1}\|_{L_2\to  L_2} \\
\le &\|\widetilde \Delta_j\mathcal E_0-  \mathcal E_0 \|_{L_2\to
L_2} + \|\widetilde \Delta_j\mathcal E_{-1}-  \mathcal E_{0} \|_{L_2\to
L_2}.
\end{align*}
The required estimate (\ref{4444}) (when $k=0$) will be a consequence   of
the pair of inequalities
\begin{align}
& \|\mathcal E_0\widetilde \Delta_j- \widetilde \Delta_j \|_{L_2\to
L_2} + \|\mathcal E_{-1}\widetilde \Delta_j- \widetilde \Delta_j \|_{L_2\to
L_2} \le C 2^{j} &\text{when $j\le 0$,}\label{c1} \\
  &\|\widetilde \Delta_j\mathcal E_0-  \mathcal E_0 \|_{L_2\to
L_2} + \|\widetilde \Delta_j\mathcal E_{-1}-  \mathcal E_{0} \|_{L_2\to
L_2}\le C 2^{-j} &\text{when $j\ge 0$.} \label{c2}
\end{align}
We start by proving (\ref{c1}). We only consider the term
$\mathcal E_0\widetilde \Delta_j- \widetilde \Delta_j$ since the
term $\mathcal E_{-1}\widetilde \Delta_j- \widetilde \Delta_j$ is similar. Let
$f\in L_2(\mathbb R^n)$. Then
\begin{align*}
&\|\mathcal E_0\widetilde \Delta_j f- \widetilde \Delta_jf \|_{L_2}^2
=   \sum_{Q\in \mathcal D_0} \|f*\phi_j -\Ave_Q (f*\phi_j) \|_{L_2(Q)}^2\\
\le & \sum_{Q\in \mathcal D_0} \int_Q\int_Q
|(f*\phi_j)(x)-(f*\phi_j)(t)|^2\, dt\, dx\\
\le &   \sum_{Q\in \mathcal D_0} \int_Q\int_Q \Big(\int_{3Q} |f(y)|
|\phi_j (x-y)|  \,dy \Big)^2 \, dt\, dx \\
&\quad\quad\,\, + \sum_{Q\in \mathcal D_0}
\int_Q\int_Q \Big(\int_{3Q} |f(y)|   |\phi_j (t-y)|  \,dy \Big)^2 \, dt\, dx  \\
&\quad\quad\,\, +\sum_{Q\in \mathcal D_0}
\int_Q\int_Q \Big(\int_{(3Q)^c} |f(y)|  2^{jn+j} |\nabla \phi(2^j(\xi_{x,t}-y))|
\,dy \Big)^2 \, dt\, dx ,
\end{align*}
where $\xi_{x,t}$ lies on the line segment between $x$ and $t$.
It is now easy to see that the sum of the last three expressions above is bounded by
$$
C   2^{2jn} \sum_{Q\in \mathcal D_0} \int_{3Q} |f(y)|^2\, dy + C_M 2^{2j}
\sum_{Q\in \mathcal D_0} \int_Q \Big(\int_{\mathbb R^n}  \frac{2^{jn}|f(y)|\,
dy}{(1+ 2^j|x-y|)^M} \Big)^2dx
$$
which is clearly controlled by $C 2^{2j} \|f\|_{L_2}^2$. This
estimate is useful when $j\le 0$.

We now turn to the proof of (\ref{c2}).   Since $\widetilde \Delta_j$ is
the difference of two $\widetilde S_j$'s, it will suffice to prove
(\ref{c2}) where $\widetilde \Delta_j$ is replaced by $\widetilde S_j$.
We only
work with the term $\widetilde S_j\mathcal E_0-  \mathcal E_0 $
since the other term can be treated similarly.
We have
\begin{align*}
&\|\widetilde S_j\mathcal E_0 f-  \mathcal E_0f\|_{L_2}^2 =
\big\|\sum_{Q\in \mathcal D_0} (\Ave_Q f) \, (\psi_j *\chi_Q-\chi_Q) 
\big\|_{L_2}^2 
\\
\le &
2\big\|\sum_{Q\in \mathcal D_0} (\Ave_Q f) \, (\psi_j *\chi_Q-\chi_Q)
\chi_{3Q}\big\|_{L_2}^2 +2\big\|\sum_{Q\in \mathcal D_0} (\Ave_Q f) \,
(\psi_j *\chi_Q )
\chi_{(3Q)^c}\big\|_{L_2}^2.
\end{align*}
Since the functions appearing inside the sum in the first term
above have supports with bounded overlap we obtain
$$
\big\|\sum_{Q\in \mathcal D_0} (\Ave_Q f) \, (\psi_j *\chi_Q-\chi_Q)
\chi_{3Q}\big\|_{L_2}^2 \le C\sum_{Q\in \mathcal D_0} (\Ave_Q |f|)^2
\|\psi_j *\chi_Q-\chi_Q\|_{L_2}^2,
$$
and the crucial observation is that
$$
\|\psi_j *\chi_Q-\chi_Q\|_{L_2} \le C 2^{-j} ,
$$
which can be easily checked using the Fourier transform. Therefore
we obtain
$$
\big\|\sum_{Q\in \mathcal D_0} (\Ave_Q f) \, (\psi_j *\chi_Q-\chi_Q)
\chi_{3Q}\big \|_{L_2}^2 \le C2^{-2j}\|f\|_{L_2}^2,
$$
and the required conclusion will be proved if we can show that
\begin{equation}\label{55}
\big \|\sum_{Q\in \mathcal D_0} (\Ave_Q f) \,
(\psi_j *\chi_Q )
\chi_{(3Q)^c}\big \|_{L_2}^2 \le C2^{-2j} \|f\|_{L_2}^2.
\end{equation}
We prove (\ref{55}) by using a purely size estimate.
Let $c_Q$ be  the center of the dyadic cube $Q$.
For  $x\notin 3Q$  we have  the  easy estimate
$$
|(\psi_j*\chi_Q)(x)|\le \frac{C_M 2^{jn}}{ (1+2^j|x-c_Q|)^{M}} \le
 \frac{C_M 2^{jn}}{ (1+2^j )^{M/2}}\frac{1}{ (1+ |x-c_Q|)^{M/2}}
$$
since both $2^j\ge 1, |x-c_Q|\ge 1$.   We now control the left hand side of
(\ref{55}) by
\begin{align*}
&2^{j(2n-M)}
\sum_{Q\in \mathcal D_0} \sum_{Q'\in \mathcal D_0} (\Ave_Q |f|)(\Ave_{Q'} |f|)
\int_{\mathbb R^n} \!\frac{C_M\,\,dx}{ (1\!+\! |x\! -\! c_Q|)^{\frac{M}2}
  (1\!+\! |x\!-\! c_{Q'}|)^{\frac{M}2}} \\
\le &2^{j(2n-M)}
\sum_{Q\in \mathcal D_0} \sum_{Q'\in \mathcal D_0} \frac{
(\Ave\limits_Q |f|)(\Ave\limits_{Q'} |f|)}{ (1+|c_Q-c_{Q'}|)^{\frac{M}4}}
\int_{\mathbb R^n} \!\frac{C_M\,\,dx}{ (1\!+\! |x\! -\! c_Q|)^{\frac{M}4}
  (1\!+\! |x\!-\! c_{Q'}|)^{\frac{M}4}}\\
\le &2^{j(2n-M)}
\sum_{Q\in \mathcal D_0} \sum_{Q'\in \mathcal D_0} \frac{
C_M'}{ (1+|c_Q-c_{Q'}|)^{\frac{M}4}} \bigg( \int_Q |f(y)|^2\, dy+
\int_{Q'} |f(y)|^2\, dy\bigg) \\
\le &C_M'' 2^{j(2n-M)} \sum_{Q\in \mathcal D_0}\int_Q |f(y)|^2\, dy
= C_M'' 2^{j(2n-M)}\|f\|_{L_2}^2.
\end{align*}
By taking $M$ large enough we obtain (\ref{55}) and thus (\ref{c2}).
\end{proof}

We have the following result relating $h(A)$ and $\widetilde h_p(A)$.

\begin{Thm}\label{equivalence} For  every $1<p<\infty$, there is a constant $C$
depending only on
$\psi$ and $p$ so that for any $c_{00}-$matrix $A$ we have
$$
\frac1C h(A)\le\widetilde h^w_p(A)\le\widetilde h_p(A)\le  Ch (A).
$$
\end{Thm}

\begin{proof}
Consider the operators $V_r,\ r\in\mathbb Z$ defined by
$$ V_r=\sum_{j\in\mathbb Z}\Delta_j\widetilde\Delta _{j+r}.$$
Then
\begin{equation*}
V_rV_r^*  =\sum_{j,k} \Delta_j\widetilde\Delta_{j+r}\widetilde\Delta
_{k+r}\Delta_k = \sum_{|j-k|\le 1} \Delta_j
\widetilde\Delta_{j+r}\widetilde\Delta_{k+r}\Delta_k.
\end{equation*}
Hence by splitting into 3 pieces and using Proposition \ref{LP-martingale}
we obtain the estimate
$$ \|V_r\|_{L_2\to L_2} \le C2^{-|r|}.$$

Next pick $q$ so that $1<q<\infty$ and
$\frac1p=\frac{\theta}q+\frac{1-\theta}2$ where $0<\theta<1.$  Let
$(\epsilon_j)_{j\in\mathbb Z}$ be a sequence of independent Bernoulli
random variables on some probability space $(\Omega,\mathbb P).$  Then
for $f\in L_q(\Omega) $ we have
$$
 V_rf =\int_{\Omega} \sum_{j\in\mathbb Z}\sum_{k\in\mathbb Z}
\epsilon_j(\omega)\epsilon_{k-r}(\omega)\Delta_j\widetilde\Delta_kf\,
d\mathbb P.
$$
Averaging now gives
$$
\|V_rf\|_{L_q}\le (\max_{\omega}\|\sum_{j\in\mathbb
Z}\epsilon_j(\omega)\Delta_j\|_{L_q\to L_q})(\max_{\omega}\|\sum_{k\in\mathbb Z}
\epsilon_{k-r}(\omega)\widetilde \Delta_k\|_{L_q\to L_q})\|f\|_{L_q}.
$$
Hence
$\|V_r\|_{L_q\to L_q}\le C$ where $C$ depends only on $q $ and $\psi$.  Similarly
$\|V_r^*\|_{L_q\to L_q}\le C.$  By interpolation we obtain $\|V_r\|_{L_p\to L_p},
\|V_r^*\|_{L_p\to L_p}\le C2^{-|r|(1-\theta)}.$

Finally let us write
\begin{align*}
 \sup_{j\in \mathbb Z}\big|\sum_{k\in \mathbb Z}a_{jk}\widetilde\Delta_kf\big|
&= \sup_{j\in
\mathbb Z} \big|\sum_{k\in \mathbb Z}a_{jk} \sum_{r\in \mathbb Z}
\Delta_{k-r}\widetilde\Delta_kf\big|\\ &\le
\sum_{r\in \mathbb Z}\sup_{j\in
\mathbb Z}\big|\sum_{k\in \mathbb Z} a_{j,k+r} \Delta_{k}
\widetilde\Delta_{k+r}f\big| .
\end{align*}
Thus by Proposition \ref{filtration},
$$ 
\|\sup_{j\in \mathbb Z}|\sum_{k\in \mathbb Z}a_{jk}
\widetilde\Delta_kf|\|_{L_p} \le
Ch(A)\sum_{r\in
\mathbb Z} \|V_rf\|_{L_p} \le Ch_p(A) \|f\|_{L_p}.
$$
This shows that $\widetilde h_p(A)\le Ch(A).$

For the converse direction we use $V_r^*$ and  Lemma \ref{translation}.
We have
\begin{equation*}
 \sup_{j\in \mathbb Z} \big|\sum_{k\in \mathbb Z}a_{jk}\Delta_kf\big|
\le \sum_{r\in \mathbb Z}\sup_{j\in \mathbb Z}\big|\sum_{k\in \mathbb Z} a_{j,k+r}
\widetilde\Delta_{k}\Delta_{k+r}f\big|
\end{equation*}
and so
$$
\big\|\sup_{j\in\mathbb Z}\big|\sum_{k\in\mathbb Z}a_{jk}\Delta_kf
\big|\big\|_{L_{p,\infty}}\le C\widetilde h_p^w(A)\sum_{r\in\mathbb Z}
\|V_r^*f\|_{L_p}
$$
which leads  to the estimate $h(A)\le C\widetilde h_p^w(A).
$\end{proof}

We next extend the definition of $\widetilde h_p(A)$ to the case when $0<p\le
1$. For such $p$'s we define
$\widetilde h_p(A)$ to be the least constant so that for
$f\in\mathcal S $ we have
\begin{equation}\label{pless1}
\|\sup_{j\in\mathbb Z}|\sum_{k\in\mathbb
Z} a_{jk}\widetilde \Delta f|\|_{L_p}\le C\|f\|_{H_p}.
\end{equation}
The space $H_p$ that appears on the right of (\ref{pless1}) when $0<p\le 1$
is the classical   real
Hardy space of Fefferman and Stein \cite{FS} and the expression
$\|\,\,\|_{H_p}$ is its quasi-norm. 

\begin{Thm}\label{hardy}
If $0<p<1$  then there is constant $C=C(p,\psi)$
so that $C^{-1}h(A)\le \widetilde h_p(A)\le Ch(A).$
\end{Thm}

\begin{proof} First we show the estimate $\widetilde h_p(A)\le C h(A).$
Using the atomic
characterization of $H_p$, \cite{coifman}, we note that it suffices to
get an
estimate for a function $f\in \mathcal S$ supported in a cube $Q$ so that
$|f(x)|\le |Q|^{-\frac1p}$ for $x\in Q$ and
$\int x^{\alpha}f(x)=0$ if $|\alpha|\le N=[n(\frac1p-1)].$   It is then
easy to see that if $x\notin 2Q$
$$
 |\sum_{k\in\mathbb Z}a_{jk}\widetilde\Delta_kf(x)|\le
C h(A)|x-c_Q|^{-n-N-1}
$$
since $|a_{jk}|\le Ch(A)$ for each $j,k$. (Here $2Q$ is the cube with twice
the length and the same center $c_Q$ as usually.)
This gives the estimate
$$
\int_{\mathbb R^n\setminus 2Q}\sup_j|\sum_{k}a_{jk}\widetilde\Delta_kf(x)|^pdx\le
C^ph(A)^p.
$$
On the other hand,
$$ \int_{2Q}
\sup_j \big|\sum_{k}a_{jk}\widetilde\Delta_kf(x)\big|^pdx\le
C|Q|^{1-\frac{p}{2}}h(A)^p\bigg(\int_Q|f(x)|^2dx\bigg)^{\frac{p}{2}}
$$
and combining with the previous estimate we obtain $\widetilde h_p(A)\le Ch(A).$

Complex interpolation gives  that  $\widetilde
h_q(A)\le \widetilde h_2(A)^{\theta}\widetilde h_p(A)^{1-\theta}$
when   $1<q<2$   and
$\frac1q=\frac{1-\theta}{p}+\frac{\theta}{2}.$  Since $\widetilde h_q(A)\ge
C^{-1}h(A)$ we deduce the estimate $\widetilde h_p(A)\ge C^{-1}h(A).$
\end{proof}

\section{Bilinear operators}\label{bilinear}
 \setcounter{equation}0

Let $\sigma$ be a bounded measurable function on $\mathbb
R^n\times\mathbb R^n$.  For $f,g\in\mathcal S(\mathbb R^n)$ we define
a bilinear operator $W_{\sigma}(f,g)$   with multiplier $\sigma$ by setting 
\begin{equation}\label{bilineardef}
W_{\sigma}(f,g)(x) =\int_{\mathbb R^n}\int_{\mathbb R^n}
\sigma (\xi ,\eta)\widehat f(\xi )\widehat
g(\eta)e^{2\pi i \langle x ,   \xi+\eta  \rangle }\, d\xi d\eta.
\end{equation}
If (\ref{bilineardef}) is satisfied we say that $\sigma$ is the
bilinear symbol (or multiplier) of $W_\sigma$.
 Now suppose $1<p_1,p_2<\infty$ and let $p_0$ be
defined by
$\frac1{p_0}=\frac1{p_1}+\frac1{p_2}.$  We say that $W_{\sigma}$ is
strongly $(p_1,p_2)-$bounded  if $W_{\sigma}$ extends to a bounded
bilinear operator from $L_{p_1}\times L_{p_2}\to  L_{p_0}$.
In this   case we
denote its norm by $\|W_{\sigma}\|_{L_{p_1}\times L_{p_2}\to L_{p_0}}$ 
(we define this be expression to be
$\infty$ if $W_{\sigma}$ is not bounded).  Similarly we say
$W_{\sigma}$ is weakly $(p_1,p_2)-$bounded if it extends to a bounded
bilinear operator from $L_{p_1}\times L_{p_2}\to L_{p_0,\infty}$ and its
norm is then denoted $\|W_{\sigma}\|_{L_{p_1}\times L_{p_2}\to L_{p_0,\infty}}.$

We extend  these definitions to the case $0< p_1,p_2<\infty$ 
 by replacing the $L_p $ spaces   by the
corresponding Hardy spaces when $0<p_j\le 1$.  In the definition below 
we set $H_p=L_p$ for $1<p<\infty$.  Given 
$0<p_1, <p_2<\infty$ and $p_0$ defined by
$\frac1{p_0}=\frac1{p_1}+\frac1{p_2},$ we say that
$W_{\sigma}$ is strongly $(p_1,p_2)-$bounded if it extends to a
bounded bilinear operator from $H_{p_1}\times H_{p_2}\to  L_{p_0}$, and we
denote its norm by $\|W_{\sigma}\|_{H_{p_1}
\times H_{p_2}\to L_{p_0}}$.  We say that
$W_{\sigma}$ is weakly $(p_1,p_2)-$bounded if it extends to a
bounded bilinear operator from $H_{p_1}\times H_{p_2}\to  L_{p_0,\infty}$, and
in this case we denote its norm by $\|W_{\sigma}\|_{H_{p_1}\times H_{p_2}\to
L_{p_0,\infty}}$.  
Now for a bounded function $\sigma$ on $\mathbb
R^n\times \mathbb R^n$ and $ 0<p_1,p_2<\infty$ we define its 
corresponding strong and weak $(p_1,p_2)$-multiplier norm  by 
$$
 \|\sigma\|_{\mathcal M_{p_1,p_2}}=\|W_{\sigma}\|_{H_{p_1}\times H_{p_2}\to
L_{p_0}}\quad\text{and}
\quad \|\sigma\|_{\mathcal M^w_{p_1,p_2}}=\|W_{\sigma}\|_{H_{p_1}\times
H_{p_2}\to L_{p_0,\infty}},
$$
where $1/p_0=1/p_1+1/p_2$. 

This definition of multiplier norm is analogous to that 
in the linear case. 
If $\upsilon\in L_{\infty}(\mathbb R^n)$,   
$\|\upsilon\|_{\mathcal M_{p}}$ 
denotes the norm of $\upsilon$ as  a multiplier from $H_p$ into
 $L_{p }$ that is 
$$
 \|\upsilon\|_{\mathcal M_{p }}= \|
M_\upsilon   \|_{H_p\to L_p }, \quad\quad\text{where}\quad 
M_\upsilon f= (\upsilon
\widehat{f}\,)
\spcheck, 
$$
when $0<p<\infty$.  
Next
we mention a few properties of multipliers.

\begin{Prop}\label{inv}
Suppose $\sigma\in L_{\infty}(\mathbb R^n\times
\mathbb R^n)$ and $0<p_1,p_2<\infty.$  Then:\newline
(i) If $\sigma'(\xi,\eta)\!=\!\sigma(\xi\!-\!\xi_0,\eta 
\!-\!\eta_0)$ for some fixed
$\xi_0,\eta_0$ then
$\|\sigma'\|_{\mathcal
M_{p_1,p_2}}=\|\sigma\|_{\mathcal M_{p_1,p_2}}.$\newline
(ii) If $L:\mathbb R^n\to\mathbb R^n$ is an invertible linear operator
and $\sigma_L(\xi,\eta)=\sigma(L\xi,L\eta)$ then
$\|\sigma_L\|_{\mathcal
M_{p_1,p_2}}=\|\sigma\|_{\mathcal M_{p_1,p_2}}.$\newline
(iii) If $\mu,\upsilon\in L_{\infty}(\mathbb R^n)$ and
$\sigma'(\xi,\eta)=\mu(\xi)\sigma(\xi,\eta)\upsilon(\eta),$ then
$$\|\sigma'\|_{\mathcal M_{p_1,p_2}}\le
\|\mu\|_{\mathcal
M_{p_1}}\|\sigma\|_{\mathcal
M_{p_1,p_2}}\|\upsilon\|_{\mathcal M_{p_2}}.$$
\end{Prop}

\begin{proof} For (i) note that $W_{\sigma'}(f,g)=e^{2\pi i\langle
x,\xi_0+\eta_0\rangle}W(e^{-2\pi i\langle x,\xi_0\rangle}f,e^{-2\pi i\langle
x,\eta_0\rangle}g).$  For (ii) note that   $W_{\sigma_L}(f,g)\!\circ\!
(L^t)^{-1}= W(f\!\circ\!
(L^t)^{-1},g\!\circ \!(L^t)^{-1}).$    (iii) is trivial.
\end{proof}



\begin{Lem}\label{bounded2}
Let $\sigma\in L_{\infty}(\mathbb R^n\times\mathbb R^n)$.   Suppose that
either $p_0\ge 1,$ or that $\sigma$ is locally
Riemann-integrable (i.e. continuous except on a set of measure zero).
Then  $\|\sigma\|_{L_\infty}\le
 \|\sigma\|_{\mathcal M_{p_1,p_2}}$ whenever   $p_0= p_1p_2/(p_1+p_2)$ 
and $0<p_1,p_2<\infty.$ 
\end{Lem}

\begin{proof}
   Suppose that $\sigma$ is locally
Riemann-integrable and let $(\xi_0,\eta_0)$ be a point
of continuity of $\sigma.$  Then  if we put
$\sigma'_{\lambda}(\xi,\eta)= \sigma
(\xi_0+\lambda\xi,\eta_0+\lambda\eta),$   Proposition \ref{inv} gives that
$\|W_{\sigma'_{\lambda}}\|_{H_{p_1}\times H_{p_2} \to
L_{p_0}} =\|W_{\sigma}\|_{H_{p_1}\times H_{p_2} \to
L_{p_0}}.$  Now if
$f,g\in\mathcal S$ it is easy to see that as $\lambda\to 0$  we have
convergence in $L_2$ (and even pointwise) of $W_{\sigma'_{\lambda}}(f,g)$
to $ \sigma(\xi_0,\eta_0)f(x)g(x).$

If $p_0\ge 1$ let $Q_k$ be a cube of side $2^{-k}$ centered at $(0,0)$
in $\mathbb R^n\times \mathbb R^n.$  Let
$$
\sigma_k(\xi,\eta)=\frac{1}{|Q_k|}\int_{Q_k}
\sigma(\xi+\xi_0,\eta+\eta_0)d\xi_0\,d\eta_0.
$$
Proposition \ref{inv} and the fact that $p_0\ge 1$ easily imply that
$\|W_{\sigma_k}\|_{L_{p_1}\times L_{p_2} \to
L_{p_0}}\le \|W_{\sigma}\|_{L_{p_1}\times L_{p_2} \to
L_{p_0}}.$  Since
$\sigma_k$ is continuous we have
$\|\sigma_k\|_{L_\infty}\le  \|W_{\sigma}\|_{L_{p_1}\times L_{p_2} \to
L_{p_0}}$.  Taking
limits as $k\to \infty$ yields the conclusion.
\end{proof}

Next we   require a lemma on series in $L_p.$

\begin{Lem}\label{easy}
Let $0<p<\infty.$  Suppose  that for some
$(f_{jk})_{(j,k)\in\mathbb Z^2 }$   sequence  of $L_p$ functions
 and  for all pairs of sequences
$(\delta_j)_{j\in\mathbb Z}, (\delta'_k)_{k\in\mathbb Z}$ with
$\sup_{j\in\mathbb Z}|\delta_j|\le 1$ and $\sup_{j\in\mathbb
Z}|\delta_j'|\le 1,$
we have
$$
\sup_{N\in \mathbb N}  \big\|\sum_{|j|\le N}\sum_{|k|\le
N}\delta_j\delta'_kf_{jk}\big\|_{L_p} \le M.
$$
Then there is a constant $C=C(p)$ such that
\newline (i) $\sup_{|\delta_j|\le 1}\|\sum_{j\in\mathbb
Z}\delta_jf_{jj}\|_{L_p}\le CM$
(and the series converges unconditionally),\newline
(ii) $\|(\sum_{j\in\mathbb Z}\sum_{k\in\mathbb
Z}|f_{jk}|^2)^{\frac12}\|_{L_p} \le CM.$
\end{Lem}

\begin{proof}  In fact (ii) follows immediately from Khintchine's
inequality by taking  $\epsilon_j,\epsilon'_k$   two
mutually independent sequences of Bernoulli random variables.  To
obtain (i), take $\epsilon_j$ be a sequence of Bernoulli
random variables and for any finite subset
$\mathcal F\subset\mathbb Z $   write
\begin{equation}\label{kkkk-ll}
 \sum_{j\in \mathcal F}\delta_jf_{jj} = \sum_{j\in\mathcal
F}\sum_{k\in\mathcal
F}\delta_j\epsilon_j\epsilon_kf_{jk}-\sum_{\substack{j,k\in\mathcal F\\
j\neq k}}\delta_j\epsilon_j\epsilon_kf_{jk}.
\end{equation}
Now for all $|\delta_j|\le 1$, (see  also \cite{K1}, proof of Theorem 4.6),
\begin{align*}
&\mathbb E(\|\sum_{\substack{j,k\in\mathcal F\\ j\neq
k}}\delta_j\epsilon_j\epsilon_kf_{jk}\|_{L_p}^p)^{1/p}\le
C\big\|(\sum_{\substack{j,k\in\mathcal F\\
j<k}} |\delta_jf_{jk}+\delta_kf_{kj}|^2)^{\frac12}\big\|_{L_p}  \\ \le
&C\big\|(\sum_{j\in\mathbb Z}\sum_{k\in\mathbb
Z}|f_{jk}|^2)^{\frac12}\big\|_{L_p}  \le CM
\end{align*}
by a generalization of Khintchine's inequality due to Bonami
\cite{bonami} and part (ii).
The same estimate is also valid for $\sum\limits_{j\in\mathcal
F}\sum\limits_{k\in\mathcal
F}\delta_j\epsilon_j\epsilon_kf_{jk}$ by our assumptions.  These
estimates together with
  (\ref{kkkk-ll})   give       (i).
\end{proof}

We now introduce some notation that will be useful in the sequel.
For $(j,k)\in\mathbb Z$ let $D_{jk}=
\{(\xi,\eta):2^{j-1}\le |\xi|\le
2^{j+1},\ 2^{k-1}\le |\eta|\le 2^{k+1}\}.$
Also for $\theta>0$
 let $D_{jk}(\theta)=
\{(\xi,\eta):2^{j-\theta}\le |\xi|\le
2^{j+\theta},\ 2^{k-\theta}\le |\eta|\le 2^{k+\theta}\}.$

\begin{Prop}\label{subdiagonal} For any $1<p_1,p_2<\infty$ there is a
constant
$C=C(p_1,p_2)$ so that whenever $(\sigma_{jk})_{j,k\in\mathbb Z}$ is a
family of bilinear symbols  with $\supp
\sigma_{jk}\subset D_{jk}$   which satisfy 
$$
\sup_{  |\delta_j|\le 1}
\sup_{  |\delta_k'|\le 1}\|\sum_{j}\sum_k\delta_j\delta'_k
\sigma_{jk}\|_{\mathcal M_{p_1,p_2}}\le M,
$$
  then the following statements are valid:\newline
(i) For any scalar sequence $(\delta_j)$
 with $\sup_j|\delta_j|\le 1$ and any $r\in\mathbb Z$ we have
$$
\|\sum_{j\in\mathbb Z}\delta_j \sigma_{j,j+r}\|_{\mathcal M_{p_1,p_2}}\le CM.
$$
(ii)  For all $r\ge 3$ we have,
$$
 \|\sum_{j\in\mathbb Z}\sum_{k\le
j-r}\sigma_{jk}\|_{\mathcal M_{p_1,p_2}}+ \|\sum_{k\in\mathbb
Z}\sum_{j\le k-r}\sigma_{jk}\|_{\mathcal M_{p_1,p_2}}\le
C(1+r^{\max(\frac1{p_0},1)})M.
$$
(iii) For every $r\ge 3$, $p_0\le 1$ and for all
$f,g\in\mathcal S,$ we have
\begin{align*}
&\big\|\sum_{j\in\mathbb Z}\sum_{k\le j-r}W_{\sigma_{jk}} \big\|_{
L_{p_1}\times L_{p_2} \to H_{p_0}} 
 \le C(1+r^{\max(\frac1{p_0},1)})M \\
&\big\|\sum_{k\in\mathbb
Z}\sum_{j\le k-r}    W_{\sigma_{jk}} \big\|_{
L_{p_1}\times L_{p_2} \to H_{p_0}}   
 \le C(1+r^{\max(\frac1{p_0},1)})M .
\end{align*}
\end{Prop}

\begin{proof}
For simplicity we write $W_{jk}=W_{\sigma_{jk}} $ below. (i) follows directly
from Lemma \ref{easy}.  To prove (ii) and (iii)
it is enough to consider the case $r=3$, since the
other cases follow trivially by applying (i) and the known case $r=3.$ We
therefore suppose $r\ge 3$ and establish both (ii) and (iii).
An easy calculation gives  that
for $f,g$ Schwartz, the function $W_{jk}(f,g)$
 has Fourier transform supported
in the annulus $ 2^{j-2}\le |\zeta |\le 2^{j+2}$ when $k\le j-3$.
It follows that
\begin{align}\begin{split}\label{kl}
\|\sum_{j\in\mathbb Z}\sum_{k\le j-3}W_{jk}(f,g)\|_{L_{p_0}} &\le
\|\sum_{j\in\mathbb Z}\sum_{k\le j-3}W_{jk}(f,g)\|_{H_{p_0}} \\
&\le C\|(\sum_{j\in\mathbb Z}|\sum_{k\le j-3}W_{jk}(f,g)|^2
)^{\frac12}\|_{L_{p_0}} \\
&\le C\mathbb E (\|\sum_{j\in\mathbb Z}\epsilon_j
\sum_{k\le j-3}W_{jk}(f,g)\|_{L_{p_0}}^{p_0})^{1/p_0}
\end{split}\end{align}
where as usual $(\epsilon_j)_{j\in\mathbb Z}$ is a sequence of
independent Bernoulli random variables. (If $p_0>1$ then
$H_{p_0}=L_{p_0}$.) We need to control the last term in (\ref{kl}).

Our hypothesis gives the estimate
\begin{equation}\label{aaaa}
 \mathbb E(\|\sum_{j\in\mathbb Z}\sum_{k\in\mathbb
Z}\epsilon_jW_{jk}(f,g)\|_{L_{p_0}}^{p_0})^{1/p_0} \le
CM\|f\|_{L_{p_1}}\|g\|_{L_{p_2}},
\end{equation}
while we can apply (i) to obtain
\begin{equation}\label{bbbb}
 \mathbb E(\|\sum_{j\in\mathbb Z}\sum_{|k-j|\le
2}\epsilon_jW_{jk}(f,g)\|_{L_{p_0}}^{p_0})^{1/p_0}\le
CM\|f\|_{L_{p_1}}\|g\|_{L_{p_2}}.
\end{equation}
It remains to estimate
\begin{align*}
 \mathbb E \big(\|\sum_{j\in\mathbb Z}\sum_{k\ge
j+3}\epsilon_jW_{jk}(f,g)\|_{L_{p_0}}^{p_0}\big)^{1/p_0}&\le
 \mathbb E\big(\|\sum_{j\in\mathbb Z}\sum_{k\ge
j+3}\epsilon_jW_{jk}(f,g)\|_{H_{p_0}}^{p_0}\big)^{1/p_0}\\ &\le
C \mathbb E\big(\|\sum_{k\in\mathbb Z}(\sum_{j\le
k-3}\epsilon_jW_{jk}(f,g)\|_{L_{p_0}}^{p_0}\big)^{1/p_0}\\
&\le C\mathbb E\big(\|\sum_{k\in\mathbb Z}\epsilon'_k\sum_{j\le
k-3}\epsilon_j W_{jk}(f,g)\|_{L_{p_0}}^{p_0}\big)^{1/p_0}
\end{align*}
where $\epsilon_k'$ is a second (independent) sequence of independent
Bernoulli random variables.
Hence using again Khintchine's inequality
we have
\begin{align}\begin{split}\label{kkkk}
 \mathbb E \big(\|\sum_{j\in\mathbb Z}\sum_{k\ge
j+3}\epsilon_jW_{jk}(f,g)\|_{L_{p_0}}^{p_0}\big)^{1/p_0}&\le
 C\|(\sum_{k\in\mathbb
Z}\sum_{j\le k-3}|W_{jk}(f,g)|^2)^{\frac12}\|_{L_{p_0}}\\
&\le C\|(\sum_{k\in\mathbb Z}\sum_{j\in \mathbb Z}
|W_{jk}(f,g)|^2)^{\frac12}\|_{L_{p_0}}\\
&\le CM\|f\|_{L_{p_1}}\|g\|_{L_{p_2}}
\end{split}\end{align}
in view of  Lemma \ref{easy}.
Using (\ref{aaaa}), (\ref{bbbb}), and (\ref{kkkk})  we obtain
$$
\mathbb E \big(\|\sum_{j\in\mathbb Z}\epsilon_j\sum_{k\le
j-3}W_{jk}(f,g)\|^{p_0}_{L_{p_0}} \big)^{1/p_0}
\le CM\|f\|_{L_{p_1}}\|g\|_{L_{p_2}}
$$
which combined with (\ref{kl}) gives  the first of the
 assertions (ii) and (iii) for $r=3$. The second assertions are derived
similarly by symmetry.
\end{proof}

We will need one further preliminary lemma.

\begin{Lem}\label{factors}
For any $1<p_1,p_2<\infty$ there is a constant
$C=C(p_1,p_2)$ such that for any family of symbols
$(\sigma_{jk})_{j,k\in\mathbb Z }$
 with
$\supp \sigma_{jk}\subset D_{jk}$ and
for any $\mu,\upsilon$
$C^{\infty} $ functions on the annulus
$\frac14\le |\xi|\le 4 $  we have
\begin{equation*}
 \sup_{|\delta_j|}\sup_{|\delta_k'|\le 1} \big\|
\sum_{j\in\mathbb Z}\sum_{k\in\mathbb Z}
\delta_j\delta_k' \tau_{jk} \big\|_{\mathcal M_{p_1,p_2}}
 \le CK_{\mu}K_{\upsilon} \sup_{|\delta_j|}\sup_{|\delta_k'|\le 1}
\big\|\sum_{j}\sum_k\delta_j\delta'_k
\sigma_{jk}\big\|_{\mathcal M_{p_1,p_2}}, 
\end{equation*}
where
$\tau_{jk}(\xi,\eta)=\mu(2^{-j}\xi)\sigma_{jk}(\xi,\eta)\upsilon(2^{-k}\eta)$,
  $$
K_{\mu}=\sup_{\substack{|\alpha|\le m\\
\frac14\le |\xi|\le 4}}\left|\frac{\partial^{\alpha}\mu}{\partial
\xi^{\alpha}}\right|,\qquad\text{ }\qquad
K_{\upsilon}=\sup_{\substack{|\alpha|\le m\\
\frac14\le |\xi|\le 4}}\left|\frac{\partial^{\alpha}\upsilon}{\partial
\xi^{\alpha}}\right|,
$$
and  $m=[(n+1)/2].$
\end{Lem}

\begin{proof}
Recalling the definition of $\phi$  from section \ref{s-harmonic} we note that the
function
$$
\big(\sum_{l=j-2}^{j+2}\widehat\phi_l(\xi) \big)
\big(\sum_{l=j-2}^{k+2}\widehat\phi_l(\eta) \big)
$$
is compactly supported and is equal to $1$ on the support of
$\sigma_{jk}(\xi , \eta)$.
For any sequence $\delta_j$ with $\sup|\delta_j|\le 1$ we
observe that
\begin{align}\label{ftre}
& \big\|\big(\sum_{j\in\mathbb
Z}\delta_j\mu(2^{-j} \xi )\big)\big(\sum_{l=j-2}^{j+2}\widehat\phi_l(\xi)\big)
\big\|_{\mathcal M_{p_1}}\le CK_{\mu}\\
& \big\|\big(\sum_{k\in\mathbb
Z}\delta_k'\mu(2^{-jk} \eta )\big)\big(\sum_{l=k-2}^{k+2}\widehat\phi_l(\eta)\big)
\big\|_{\mathcal M_{p_2}}\le CK_{\upsilon}\label{ftre2}
\end{align}
by the H\"ormander multiplier theorem.     Let
$ U_{j_1,j_2,k_1,k_2}$ be the bilinear operator with symbol
$$
\Big(\delta_{j_1} \mu(2^{-j_1} \xi)\sum_{l=j_1-2}^{j_1+2}
\widehat\phi_l(\xi)\Big)
\sigma_{j_2,k_2}(\xi,\eta)   \Big(\delta'_{k_1}
\upsilon(2^{-k_1}\eta)\sum_{l=k_1-2}^{k_1+2}\widehat\phi_l(\eta)\Big),
$$
for some fixed $|\delta_j|, |\delta_k'|\le 1$.  Let
$$
M= \sup_{|\delta_j|}\sup_{|\delta_k'|\le 1}
\big\|\sum_{j}\sum_k\delta_j\delta'_k
\sigma_{jk}\big\|_{\mathcal M_{p_1,p_2}}
$$
and let
$(\epsilon_j),(\epsilon_k')$ be two sequences of mutually independent
Bernoulli random variables.
Then for $f,g\in\mathcal S$ we have
\begin{align*}
 &\mathbb E \big(\|\sum_{j_1\in\mathbb Z}\sum_{j_2\in\mathbb
Z}\sum_{k_1\in\mathbb Z}\sum_{k_2\in\mathbb
Z}\epsilon_{j_1}\epsilon_{j_2}\epsilon_{k_1}'\epsilon_{k_2}'   
U_{j_1,j_2,k_1,k_2} (f,g)\|_{L_{p_0}}^{p_0} \big)^{\frac{1}{ p_0}} \\
&\quad\quad\quad\quad\quad\quad\quad
\quad\quad\quad\quad\quad\quad\quad
\le CMK_{\mu}K_{\upsilon}\|f\|_{L_{p_1}}\|f\|_{L_{p_2}}
\end{align*}
by our hypothesis, (\ref{ftre}), and  (\ref{ftre2}).
We now use Lemma \ref{easy} twice   to deduce that
$$
 \|\sum_{j\in\mathbb Z}\sum_{k\in\mathbb Z}U_{j,j,k,k}(f,g)\|_{L_{p_0}}\le
CK_{\mu}K_{\upsilon}M \|f\|_{L_{p_1}}\|g\|_{L_{p_2}}.
$$
This proves the required assertion.
\end{proof}

\section{Bilinear operators and infinite matrices}\label{s-bilinear+infmatrices}
\setcounter{equation}0

Recall from   section \ref{s-harmonic} that  $\phi_j(x)=2^{nj}\phi(2^jx)$ are smooth
bumps  whose Fourier transforms are supported in the annuli
$2^{j-1}\le |\xi |\le 2^{j+1}$.
In this section we will consider   symbols $\sigma$ of the form
\begin{equation}\label{elemsymbols}
 \sigma_A(\xi,\eta)=\sum_{j\in\mathbb Z}\sum_{k\in\mathbb
Z}a_{jk}\widehat{\phi_j}(\xi)\widehat{\phi_k}(\eta)
\end{equation}
 where $A=(a_{jk})_{(j,k)\in\mathbb Z^2}$ is a bounded infinite matrix.
We let $W_A=W_{\sigma_A}$ and $\|A\|_{\infty}=\sup_{j,k}|a_{jk}|.$

If $A$ is such an infinite matrix we define $A_L$ to be its
lower-triangle and $A_U$ to be its upper-triangle i.e.
$A_L=(a_{jk}\theta_{jk})_{j,k}$ and $A_U=(a_{jk}\theta_{kj})_{j,k}$ where
$\theta_{jk}=1$ if $k<j$ and $0$ otherwise.  We let $A_D$ be the diagonal
$A-A_U-A_L.$  Now define
\begin{equation}\label{defH}
 H(A)=h(A_L)+h(A_U^t)+\|A\|_{\infty} \end{equation}
Notice that $H(A)\ge \|A\|_{\infty}$ and that $H$ is a norm on
the space of $\{A:H(A)<\infty\}$ which makes it a Banach space.

Our objective will be to show that for any choice of $0<p_1,p_2<\infty$ we
have
$\|W_A\|_{H_{p_1}\times H_{p_2}\to L_{p_0}}\approx H(A).$
This will provide  us with an equivalent expression for the norm of the
multiplier $\sigma_A$ defined in (\ref{elemsymbols}).

We start by proving the simple upper estimate below.

\begin{Lem}\label{upper}
 If $0<p_1,p_2<\infty$ there is a constant
$C=C(p_1,p_2)$ so
that for any matrix $A$ we have $\|\sigma_A\|_{\mathcal M_{p_1,p_2}}\le CH(A).$
\end{Lem}

\begin{proof}
 We give the proof in the case $p_1,p_2>1$; the only real
alteration for the other cases would be to replace the appropriate
$L_{p_j}-$norm with the $H_{p_j}-$ norm and use Theorem \ref{hardy}.  Suppose
$f,g\in\mathcal S$ and consider
\begin{align}\begin{split}\label{no0}
 W_{A}(f,g)=&\sum_{j\in\mathbb
Z}\sum_{k\le j-3}a_{jk}\widetilde\Delta_jf\widetilde\Delta_kg
+\sum_{k\in\mathbb Z}\sum_{j\le
k-3}a_{jk}\widetilde\Delta_jf\widetilde\Delta_kg  \\ &\quad\quad\quad
+\sum_{j\in\mathbb Z}\sum_{k=j-2}^{j+2}
a_{jk}\widetilde\Delta_jf\widetilde\Delta_kg.
\end{split}\end{align}
We estimate the first term by noticing that for fixed $j$ the Fourier
transform of $\widetilde\Delta_jf \sum_{k\le j-3}a_{jk}\widetilde\Delta_kg$ is
contained in the set $\{\zeta:\,\, 2^{j-2}\le |\zeta |\le 2^{j+2}\}.$
Hence if $p_0>1$ we have
$$
 \|\sum_{j\in\mathbb Z}\widetilde\Delta_jf\sum_{k\le
j-3}a_{jk}\widetilde\Delta_kg\|_{L_{p_0}} \le C\|(\sum_{j\in\mathbb
Z}|\widetilde\Delta_jf|^2)^{\frac12} (|\sum_{k\le
j-3}a_{jk}\widetilde\Delta_k|^2)^{\frac12}\|_{L_{p_0}}.
$$
If $0<p_0\le 1$ we obtain the same estimate by noticing that
$$
 \|\sum_{j\in\mathbb Z}\widetilde\Delta_jf\sum_{k\le
j-3}a_{jk}\widetilde\Delta_kg\|_{L_{p_0}} \le
 \|\sum_{j\in\mathbb Z}\widetilde\Delta_jf\sum_{k\le
j-3}a_{jk}\widetilde\Delta_kg\|_{H_{p_0}}
$$
and using the corresponding square-function estimates in $H_{p_0}$. Now we
have
\begin{align}\begin{split}\label{no1}
& \big\|(\sum_{j\in\mathbb Z}|\widetilde\Delta_jf|^2|  \sum_{k\le
j-3}a_{jk}\widetilde\Delta_kg|^2)^{\frac12}\big\|_{L_{p_0}} \\
&\quad\quad\quad\quad\quad\quad\quad\quad \le   \big\|(\sum_{j\in
\mathbb Z}|\widetilde\Delta_jf|^2)^{\frac12}\sup_{j\in\mathbb Z}|\!\!
\sum_{k\le j-3}a_{jk}\widetilde\Delta_kg|\big\|_{L_{p_0}}.
\end{split}\end{align}
If we let $A_{LL}$ be the matrix with entries $a_{jk}$ if $k\le j-3$
and $0$ otherwise, then $h(A_{LL})\le h(A_L)+h(B)$ where $B$ is the
matrix with entries $a_{jk}$ if $j-2\le k\le j-1$ and $0$ 
otherwise. It is trivial
to see that one has the estimate $h(B)\le 2\|A\|_{\infty}
$ so that
$h(A_{LL})\le Ch(A_L).$  Hence (\ref{no1}) and Theorem
\ref{equivalence} give
\begin{align*}
\| \sum_{j\in\mathbb Z}\sum_{k\le
j-3}a_{jk}\widetilde\Delta_jf\widetilde\Delta_kg\|_{L_{p_0}} 
&\le C\|(\sum_{j\in\mathbb
Z}\widetilde\Delta_jf)^{\frac12}\|_{L_{p_1}}\|\sup_{j\in\mathbb Z}|
\sum_{k\in\mathbb
Z}a_{jk}\widetilde\Delta_kg|\|_{L_{p_2}}\\
&\le Ch(A_L)\|f\|_{L_{p_1}}\|g\|_{L_{p_2}}.
\end{align*}

The same argument shows that the third term in (\ref{no0}) is controlled
by $Ch(A_U^t)\|f\|_{L_{p_1}}\|g\|_{L_{p_2}}.$  The middle term in (\ref{no0}) is
easy. For $-2\le r\le 2 $ we have
 \begin{align*}
&\big\|\sum_{j\in\mathbb
Z}a_{j,j+r}\widetilde\Delta_jf\widetilde\Delta_{j+r}g\big\|_{L_{p_0}} \\
\le & \big\|(\sum_{j\in\mathbb Z}|a_{j,j+r}||\widetilde\Delta_j
f|^2)^{\frac12}\big\|_{L_{p_1}}
\big\|(\sum_{k\in\mathbb Z}|a_{j,j+r}||\widetilde\Delta_{j+r} 
g|^2)^{\frac12}\big\|_{L_{p_2}}
\\  \le  &
C\max_{j}|a_{j,j+r}|\|f\|_{L_{p_1}}\|g\|_{L_{p_2}}.
\end{align*}
Combining we obtain the required upper estimate:
$\|\sigma_A\|_{\mathcal M_{p_1,p_2}}\le CH(A).$ \end{proof}

To obtain the converse is somewhat more complicated.  First we prove a
general result which we will use in other situations as well.

\begin{Prop}\label{averaging}
 For any $1<p_1,p_2<\infty$ with $p_0=(1/p_1+1/p_2)^{-1}\!\ge\! 1$, there is a
constant
$C=C(p_1,p_2)$ with the following property. Whenever
$(\sigma_{jk})_{(j,k)\in\mathbb Z^2}$
is a family of symbols with
$\supp \sigma_{jk}\subset D_{jk}$   which satisfy 
$$
\sup_{|\delta_j|\le 1} \sup_{|\delta'_k|\le 1}
\|\sum_{j}\sum_k\delta_j\delta'_k W_{\sigma_{jk}}\|_{L_{p_1}
\times L_{p_2}\to L_{p_0}}\le M,
$$
then
$$
 \|\sigma_A\|_{\mathcal M_{p_1,p_2}} \le CM,
$$
where  $A=(a_{jk})_{j,k}$ and
$$
 a_{jk} =\int_{\mathbb R^n}\int_{\mathbb
R^n}\sigma_{jk}(2^j\xi,
2^k\eta)d\xi\,d\eta .
$$
\end{Prop}

\begin{proof} As before we write $W_{jk}=W_{\sigma_{jk}}.$  Let us
consider first the case when
$\sigma_{jk}=0$ unless
$k\le j-5.$  Let $\upsilon$ be a $C^{\infty}-$function on $\mathbb R^n$
supported on $2^{-4}\le |\xi|\le 2^4$ and such that $\upsilon(\xi)=1$ on
$2^{-3}\le |\xi|\le 2^3.$ Fix $\xi_0\in\mathbb R^n$ and consider
the symbol
$$
\tau_{jk}(\xi_0;\xi,\eta)=\upsilon
(2^{-j}\xi)\sigma_{jk}(\xi+2^j\xi_0,\eta).
$$
Note that  $\tau_{jk}$ is supported in $D_{jk}(4).$
Let $T_{jk}$ be  bilinear operator with symbol $\tau_{jk}$.
For any sequences $(\delta_j)_{j\in\mathbb Z},(\delta'_k)_{k\in\mathbb
Z}$ with $\sup|\delta_j|,\sup|\delta'_k|\le 1$ and $f,g\in\mathcal S $ we have
$$
 \|\sum_{j\in\mathbb Z}\sum_{k\in\mathbb Z}\delta_j\delta_k'
T_{jk}(f,g)\|_{L_{p_0}} \le C\|(\sum_{j\in\mathbb Z}|\sum_{k\in\mathbb
Z}\delta_k'T_{jk}(f,g)|^2)^{\frac12}\|_{L_{p_0}}
$$
by considering the supports of the Fourier transforms.  But then for
fixed $j,$
$$
 \sum_{k\in\mathbb Z}\delta_k'T_{jk}(f,g)(x)=e^{- 2\pi i \langle
x,2^j\xi_0\rangle}\sum_{k\in\mathbb Z}\delta_kW_{jk}(f,g)(x),
$$
hence
\begin{align*}
 \|\sum_{j\in\mathbb Z}\sum_{k\in\mathbb Z}\delta_j\delta_k'
T_{jk}(f,g)\|_{L_{p_0}} &\le C\|(\sum_{j\in\mathbb Z}|\sum_{k\in\mathbb
Z}\delta_k'W_{jk}(f,g)|^2)^{\frac12}\|_{L_{p_0}}\\
& \le C \|\sum_{j\in\mathbb Z}\sum_{k\in\mathbb
Z}\delta_k'W_{jk}(f,g)\|_{H_{p_0}}\\
&\le CM\|f\|_{L_{p_1}}\|g\|_{L_{p_2}}
\end{align*}
using Proposition \ref{subdiagonal}.

Now note that if $|\xi_0|> 18$ then all $T_{jk}$ vanish.  Since $p_0\ge
1$, we integrate over $|\xi_0|\le 18$ to obtain symbols
$$
\tau'_{jk}(\xi,\eta)=\int_{|\xi_0|\le 18}
\tau_{jk}(\xi,\eta)\,d\xi_0= \upsilon(2^{-j}\xi)\int_{\mathbb
R^n}\sigma_{jk}(\xi+2^j\xi_0,\eta)d\xi_0
$$
with corresponding bilinear
operators $T_{jk}'$  satisfying
$$
 \|\sum_{j\in\mathbb Z}\sum_{k\in\mathbb
Z}\delta_j\delta'_kT'_{jk}\|_{L_{p_1}
\times L_{p_2}\to L_{p_0}}\le CM
$$
whenever $|\delta_j|,|\delta_k'|\le 1.$

Note that $\tau'_{jk}$ is supported on $D_{jk}(3).$  Also if $2^{j-3}\le
|\xi|\le 2^{j+3}$ we have that $\tau_{jk}'(\xi,\eta)$ is constant in
$\xi.$

Next let $O_n$ be the orthogonal group of $\mathbb R^n$ and let
$dL$ denote the Haar measure on this group.  Define
$$
 \tau^{\#}_{jk}(\xi,\eta)=
\int_{\frac14}^{4}\lambda^{n-1}\int_{O_n}\tau'_{jk}(\lambda L\xi,\lambda
L\eta)dL\, d\lambda ,
$$
and let $T^{\#}_{jk}$ be  the corresponding bilinear operator. If $(\xi,\eta)\in
D_{jk}$ we can compute that
$$
 \tau^{\#}_{jk}(\xi,\eta)= c2^{nk}|\eta|^{-n}a_{jk}
$$
where $c$ is a
constant depending only on dimension.
On the other hand, since $p_0\ge 1$,   Proposition \ref{inv} (ii) gives  that
$$
\|\sum_{j\in\mathbb Z}\sum_{k\in\mathbb
Z}\delta_j\delta'_kT^{\#}_{jk}\|_{L_{p_1}
\times L_{p_2}\to L_{p_0}}\le CM
$$
whenever $|\delta_j|,|\delta_k'|\le 1.$

Note that $\supp\tau^{\#}_{jk}\subset D_{jk}(6).$
Let us take $\mathbb M_1$ and $\mathbb M_2$ to be residue classes modulo
10.  Then if we replace $\delta_j$ by $\delta_j\chi_{\mathbb M_1}(j)$ and
$\delta'_k$ by $\delta'_k\chi_{\mathbb M_2}(k)$ we obtain a bilinear operator
whose symbol coincides with $a_{jk}2^{nk}|\eta|^{-n}\delta_j\delta_k'$ on
$D_{jk}$ for
$(j,k)\in\mathbb M_1\times\mathbb M_2.$  Using Proposition
\ref{inv} (iii) and the multipliers $\sum_{j\in\mathbb M_1}\widehat{\phi_j}$
and $\sum_{k\in\mathbb M_2}\widehat{\phi_k}$ we obtain that the
bilinear operator $V$ with symbol
$$
\sum_{j\in\mathbb M_1}\sum_{k\in\mathbb
M_2}\delta_j\delta_k'2^{nk}|\eta|^{-n}a_{jk}\widehat{\phi_j}(\xi)\widehat
\phi_k(\eta),
$$
satisfies $\|V\|_{L_{p_1}
\times L_{p_2}\to L_{p_0}}\le CM.$
Summing over 100 different pairs of residue classes gives a similar
estimate for the symbol
$$
\sum_{j\in\mathbb Z}\sum_{k\in\mathbb
Z}\delta_j\delta_k'2^{nk}|\eta|^{-n}a_{jk}\widehat{\phi_j}(\xi)\widehat{
\phi_k}(\eta).
$$
The last step is to remove the factor
$2^{nk}|\eta|^{-n}.$  But this can be done by using Lemma \ref{factors}
since $|\eta|^{-n}$ is $C^{\infty}$ on $\frac14\le |\eta|\le 4.$
\end{proof}

We will use this result to make an important estimate on the effect of
translation in the computation of $\|W_A\|_{L_{p_1}
\times L_{p_2}\to L_{p_0}}.$ Let us define
$A^{[r,s]}$ to be the matrix $(a_{j+r,k+s})_{j,k}.$

\begin{Lem}\label{translation2}(i)
There is a constant $C$ so that for all matrices $A$ we have
$$
 \|\sigma_{A^{[r,s]}}\|_{\mathcal M_{2,2}}\le
C^{|r-s|}\|\sigma_A\|_{\mathcal M_{2,2}}
$$
\newline
(ii)  For all $1<p_1,p_2<\infty$ with $p_0=p_1p_2/(p_1+p_2)\ge 1$, 
there is a constant $C=C(p_1,p_2)$ so that if $|\delta_j|,|\delta'_k|\le 1$ then
$$
\|\sum_{j\in\mathbb Z}\sum_{k\in\mathbb
Z}\delta_j\delta_k'a_{jk}\widehat\phi(2^{-j}\xi)
\widehat\phi(2^{-k}\eta)\|_{\mathcal
M_{p_1,p_2}}\le C \|\sigma_A\|_{\mathcal M_{p_1,p_2}},
$$  i.e.
$\|\sigma_D\|_{M_{p_1,p_2}}\le C \|\sigma_A\|_{\mathcal M_{p_1,p_2}},$
where $D=(d_{jk})_{j,k}=(\delta_j\delta_k'a_{jk})_{j,k}.$
\end{Lem}

\begin{proof}
It is clear from Proposition \ref{inv} that for
any $r\in\mathbb Z$ we have 
$$
\|W_{A^{[r,r]}}\|_{L_2\times  L_2\to L_1}=\|W_A\|_{L_2\times  L_2\to L_1 }.
$$
Thus it suffices to consider the case $r=0$ and $s=\pm 1$ and establish
a bound in this case.  To do this we consider the symbols
$$
\sigma_{jk}(\xi,\eta)=\sigma_A(\xi,\eta)\mu(2^{-j}\xi)\upsilon(2^{-k}\eta)
\widehat{\phi_j}(\xi)\widehat{\phi_k}(\eta),
$$
where $\mu,\upsilon$ are
$C^{\infty}-$functions satisfying $|\mu(\xi)|,|\upsilon(\eta)|\le 1$ for
all $\xi,\eta.$
Since $\|\sum_{j\in\mathbb
Z}\delta_j\mu(2^{-j}\xi)\widehat{\phi_j}(\xi)\|_{\mathcal M_{2}}$ is bounded by
$3$ whenever $\sup_j|\delta_j|\le 1,$ and there is a similar
bound
for $\sum_{k\in\mathbb Z}\delta'_k\upsilon(2^{-k}\eta)\widehat{\phi_k}(\eta)$
we have an immediate estimate;
$$
 \|\sum_{j\in\mathbb Z}\sum_{k\in\mathbb
Z}\delta_j\delta'_kW_{\sigma_{jk}}\|_{L_2\times L_2\to L_1} \le
9\|W_A\|_{L_2\times L_2\to L_1}.
$$
Now let
$$
 b_{jk}= \int_{\mathbb R^n}\int_{\mathbb
R^n}\sigma_{jk}(2^j\xi,2^k\eta)d\xi\,d\eta.
$$
Then we can compute
$$ b_{jk} = \sum_{r=-1}^1\sum_{s=-1}^1 c_{rs}a_{j+r,k+s}$$
where
$$
 c_{rs}= \int_{\mathbb R^n}\int_{\mathbb R^n}\mu(\xi)\upsilon(\eta)
\widehat{\phi_{-r}}(\xi)\widehat{\phi_{-s}}(\eta)\widehat{\phi_0}(\xi)
\widehat{\phi_0}(\eta)
d\xi\,d\eta.
$$

Since the functions $\widehat{\phi_{r}}$ for $-1\le r\le 1$ are linearly
independent on the support of $\widehat{\phi_0}$ we can use the above
estimate for a linear combination of a finite number of choices of
$\upsilon$ and
$\xi$ so that $c_{rs}=0$ except when $r=0$ and $s=1$, so that
$B=cA^{[0,1]}$  for some fixed constant $c\neq 0.$  By Proposition
\ref{averaging} we have $\|W_B\|_{L_2\times L_2\to L_1}\le C\|W_A\|_{L_2\times
L_2\to L_1}.$  This and the similar argument for the case $s=-1$ gives the
result (i).

For (ii) we observe that the above argument actually also yields a bound
on
$\|W_D\|_{L_2\times L_2\to L_1}$ 
when $D=(d_{jk})=(\delta_j\delta_k'b_{jk})$ (since
$\delta_j\delta'_k\sigma_{jk}$ also verifies the hypotheses of
Proposition \ref{averaging}.  By choosing a similar linear combination
we can then ensure that $b_{jk}=ca_{jk}$ and obtain the desired result.
\end{proof}

The next step is to consider a discrete model of the bilinear
operator
$W_{\sigma_A}.$  We restrict ourselves to $p_1=p_2=2$ for this, although
our calculations can be done in more generality.  If $A$ is a
$c_{00}-$matrix   we define $V_{A}:\,\, L_2\times L_2\to L_1$ by
$$
 V_{A}(f,g)=\sum_{j\in\mathbb Z}\sum_{k\in\mathbb
Z}a_{jk}\Delta_jf\Delta_kg, 
$$
where $\Delta_j$ are the martingale difference operators as 
defined in section 4. 
We then have
\begin{Lem}\label{discrete}
There is a constant $C$ so that if $A$ is a
(strictly) lower-triangular matrix we have
$ h(A)\le C\|V_{A}\|_{L_2\times L_2\to L_1}.$
\end{Lem}

\begin{proof}
This is  a stopping time argument.  Suppose $f\in L_2$
with $\|f\|_{L_2}=1.$  Note that for each $j$ the function
$f_j=\sum_{k\in\mathbb Z}a_{jk}\Delta_kf$ is $\Sigma_{j-1}$-measurable
where $\Sigma_{j-1}$ is the $\sigma-$algebra generated by the dyadic
cubes in $\mathcal D_{j-1}.$  Fix $\lambda>0$.  For each $j$ let
$\mathcal Q_j$ be the collection of cubes $Q\in\mathcal D_{j-1}$ so that
$|f_j|>\lambda$ on $Q$ and for each $j_1<j$ we have $|f_{j_1}|\le
\lambda$ on $Q.$ It is not difficult to see that
$$
 \{x:\ \max_{j\in\mathbb Z}|f_j(x)|>\lambda\}= \bigcup\limits_{j\in\mathbb
Z}\bigcup\limits_{Q\in\mathcal Q_j}Q
$$ and this is a disjoint union.  Also note the
left-hand side has finite measure.

For each $j$ be $u_j$ be a $\Sigma_j-$measurable function such that
$|u_j|=1$ everywhere and $\mathcal E_{j-1}u_j=0.$  Let
$$
 g= \sum_{j\in\mathbb Z}u_j\sum_{Q\in\mathcal Q_j}\chi_Q.
$$
Then
$$
 \|g\|_{L_2}^2= |\{x:\max_{j\in\mathbb Z}|f_j(x)|\}|
$$
and
\begin{equation*}
V_A(f,g)  = \sum_{j\in\mathbb Z}f_j\Delta_jg
 = \sum_{j\in\mathbb Z}f_ju_j\sum_{Q\in\mathcal Q_j}\chi_Q.
\end{equation*}
Hence
$$
 |V_A(f,g)|\ge \lambda \chi_{(\max_j|f_j|>\lambda)}
$$
so that we have
$$
 \lambda|\{\max_j|f_j|>\lambda\}\le  \|V_A\|_{L_2\times L_2\to L_1}.
$$
This implies that $h_2^w(A)\le \|V_A\|_{L_2\times L_2\to L_1}$ and the result
follows from Theorem \ref{equiv}.
\end{proof}

We are now ready for the main result:

\begin{Thm}\label{matrixnorms}
  Suppose $0<p_1,p_2<\infty$.  Then there
is a constant $C=C(p_1,p_2)$ so that for any infinite matrix $A$ we have
$$
 \frac1CH(A)\le \|\sigma_A\|_{\mathcal M_{p_1,p_2}^w}\le
\|\sigma_A\|_{\mathcal M_{p_1,p_2}}\le CH(A).
$$
\end{Thm}

\begin{proof} The upper bound is proved in Lemma \ref{upper} so we only
need to prove the lower bound. It suffices to prove the results for the
case when $A$ is a $c_{00}-$matrix. We start by considering the case
$p_1=p_2=2,$ when $A$ is strictly lower-triangular.

In this case let us estimate the norm of the discrete model $V_A.$  In
fact
\begin{align*}
V_A(f,g) &= \sum_{j\in\mathbb
Z}\sum_{k\in\mathbb Z}a_{jk}\Delta_jf\Delta_kg \\
&= \sum_{r\in\mathbb Z}\sum_{s\in\mathbb Z}\sum_{j\in\mathbb
Z}\sum_{k\in\mathbb Z}a_{jk}\widetilde\Delta_{j-r}\Delta_jf
\widetilde\Delta_{k-s}\Delta_kg\\ &= \sum_{r\in\mathbb Z}\sum_{s\in\mathbb
Z}a_{j+r,k+s}\widetilde\Delta_j\Delta_{j+r}f\widetilde\Delta_k\Delta_{k+s}g\\
&= \sum_{r\in\mathbb Z}\sum_{s\in\mathbb
Z}W_{A^{[r,s]}}(\sum_{j\in\mathbb Z}
\widetilde\Delta_j\Delta_{j+r}f,\sum_{k\in\mathbb
Z}\widetilde\Delta_k\Delta_{k+s}g)\\
&= \sum_{r\in\mathbb Z}\sum_{s\in\mathbb
Z}W_{A^{[r,s]}}(V_{-r}^*f,V_{-s}^*g),
\end{align*}
where $V_r$ is defined in
 the proof of Theorem \ref{equivalence}. Using Proposition \ref{LP-martingale}
we obtain
$$
 \|V_A\|_{L_2\times L_2\to L_1}\le C
\sum_{r\in\mathbb Z}\sum_{s\in\mathbb
Z}2^{-|r|-|s|}\|W_{A^{[r,s]}}\|_{L_2\times L_2\to L_1}.
$$
(All these  quantities are finite since $A$ has only finitely many
non-zero entries, and so there is a uniform bound on $W_{A^{[r,s]}}.$)

It follows that we have an estimate (for a suitable $C_0,$)
\begin{equation}\label{est1}
h(A)\le C_0 \sum_{r\in\mathbb Z}\sum_{s\in\mathbb
Z} 2^{-|r|-|s|}\|W_{A^{[r,s]}}\|_{L_2\times L_2\to L_1}.
\end{equation}

Next we estimate $H(A^{[r,s]}).$  If $s\ge r$ it is clear that $A$
remains lower-triangular and the invariance properties of $h(A)$ imply
that $H(A^{[r,s]})\le H(A).$  If $s<r$ then it is easy to estimate
$$
h(A^{[r,s]}_L) \le h(A_L)+(r-s)\|A\|_{\infty}
$$
and
$$
h((A^{[r,s]}_U)^t) \le (r-s)\|A\|_{\infty}.
$$
We deduce that
$$
H(A^{[r,s]})\le h(A)+|r-s|\|A\|_{\infty}
$$
for all $r,s.$
Thus we have for a suitable constant $C_0$
\begin{equation}\label{est2}
\|W_{A^{[r,s]}}\|_{L_2\times L_2\to L_1}\le C_1(1+|r-s|)h(A).
\end{equation}

Now we may pick an integer $N$ large enough so that
$$
C_1C_0\sum_{|r|>N}\sum_{|s|>N}(1+|r-s|)2^{-|r|-|s|}\le \frac12.
$$
Then we can combine (\ref{est1}) and (\ref{est2}) to obtain
\begin{equation}\label{est3}
h(A) \le C_2\sum_{|r|\le N}\sum_{|s|\le
N}\|W_{A^{[r,s]}}\|_{L_2\times L_2\to L_1}.
\end{equation}
At this point Lemma \ref{translation2} gives   the conclusion that
$$
h(A)\le C\|W_A\|_{L_2\times L_2\to L_1}.
$$
Now suppose $A$ is arbitrary.  If we let $W_{jk}$ be the bilinear operator
with symbol $a_{jk}\widehat{\phi_j}(\xi)\widehat{\phi_k}(\eta),$ Lemma
\ref{translation2} (ii) implies that we can use Proposition
\ref{subdiagonal} (ii) to deduce that  $\|W_{A_L}\|_{L_2\times L_2\to L_1}\le
C\|W_A\|_{L_2\times L_2\to L_1}$ for some absolute constant $C.$  Thus the above
argument yields $h(A_L)\le C\|W_A\|_{L_2\times L_2\to L_1}.$  Similarly
$h(A_U^t)\le C\|W_A\|_{L_2\times L_2\to L_1}$ and Lemma \ref{bounded2} is enough
to show that
$\|A\|_{\infty}\le C\|W_A\|_{L_2\times L_2\to L_1}.$  Combining these we have the
estimate
$$
H(A)\le C\|W_A\|_{L_2\times L_2\to L_1}.
$$

The proof is completed by a simple interpolation technique.  We will
argue first that an estimate of the type
\begin{equation}\label{teydr}
 H(A)\le C(p_1,p_2)\|\sigma_A\|_{\mathcal M_{p_1,p_2}}
\end{equation}
for some fixed
$1<p_1,p_2<\infty$ implies the estimate
\begin{equation}\label{huy}
H(A)\le C(q,p_2)\|\sigma_A\|_{\mathcal M_{p_1,q}^w }
\end{equation}
for every
$1<q<\infty.$  We only need to  consider the first case and $q\neq p_2$ (when
$q=p_2$ one repeats the step).
Then we may find $1<r<\infty$ and $0<\theta<1$ so that
$$
 \frac1p_2=\frac{1-\theta}{q}+\frac{\theta}{r}.
$$
The Marcinkiewicz interpolation theorem yields
\begin{equation}\label{765rt}
 \|\sigma_A\|_{\mathcal M_{p_1,p_2}}\le C(p_1,p_2,\theta)
(\|\sigma_A\|_{\mathcal M_{p_1,q}^w })^{1-\theta}(\|\sigma_A\|_{
\mathcal M_{p_2,r}  })^{\theta} .
\end{equation}
Since $\|\sigma_A\|_{\mathcal M_{p_2,r}}\le C(p_2,r)H(A)$, 
using (\ref{765rt}), and (\ref{teydr}) we obtain   estimate  (\ref{huy}) as
required (recall that we assume $A$ is a
$c_{00}$-matrix so that all these quantities are finite).

Repeated use of this argument starting from $p_1=p_2=2$ gives the theorem
in the cases $1<p_1,p_2<\infty$.

Finally in the case where either $p_1\le 1$ or $p_2\le 1$ (or both) one
can use complex interpolation to deduce
$$
\|\sigma_A\|_{\mathcal M^w_{q_1,q_2}}\le C
(\|\sigma_A\|_{\mathcal M^w_{p_1,p_2}})^{1-\theta}(\|\sigma_A\|_{
\mathcal M_{2,2}})^{\theta}
$$
where $q_1,q_2>1$ and
$$
 \frac1q_1=\frac{1-\theta}{p_1}+\frac\theta2,\qquad
\frac1q_2=\frac{1-\theta}{p_2} +\frac\theta2.
$$
This clearly extends the lower estimate to the cases $p_1,p_2\le 1.$
\end{proof}

\section{Applications to  bilinear multipliers}\label{s-applications}
\setcounter{equation}{0}

We will now consider the boundedness of the bilinear operator $W_{\sigma}$
under conditions of Marcinkiewicz type on the symbol $\sigma$.
We will say that a symbol
$\sigma$ is $C^N$ if it is $C^N$ on the set $\{(\xi,\eta):\
|\xi|,|\eta|>0\}.$    We first give an example to show that
conditions (\ref{hypothesis}) for a function $\sigma$
on $\mathbb R^{2n}$ do not imply boundedness for
the corresponding bilinear map  on
$\mathbb R^{n} \times \mathbb R^{n}$.

\noindent {\bf Example.}  There is a $C^{\infty}-$symbol $\sigma$ so that for every
pair of multi-indices $(\alpha,\beta)$ there is a constant
$C_{\alpha,\beta}$  so that
\begin{equation}\label{Marcinkiewicz}
|\xi|^{|\alpha|}|\eta|^{|\beta|}|\partial_{\xi}^{\alpha}
\partial_{\eta}^{\beta} \sigma (\xi,\eta)|\le
C_{\alpha,\beta}\end{equation}
 but $W_{\sigma}$ is not of weak type
$(p_1,p_2)$ for any $0<p_1,p_2<\infty.$

Indeed if we let $A$ be  a bounded infinite matrix and
$\sigma(\xi,\eta)=\sigma_A(\xi,\eta)$, then $\sigma$ satisfies the
condition (\ref{Marcinkiewicz}).  However $W_A$ is of weak type
$(p_1,p_2)$ if and only if $H(A)<\infty$ by theorem \ref{matrixnorms}.
At the end of Section \ref{estimates} we showed that there are
examples (with $A$ lower-triangular) where $H(A)=\infty.$

In fact more is true.  It is shown that
the condition $0<\theta<\frac12$
in (\ref{logcondition})
is insufficient to give a bound on $h(A)$ or $H(A)$
when $A$ is lower-triangular.  This means that if $0<\theta<\frac12$ we
can construct a symbol
$\sigma$ which is $C^{\infty}$, with $W_{\sigma}$ not of weak type
$(p_1,p_2)$ for any $0<p_1,p_2<\infty$ and such that for each pair of
multi-indices $(\alpha,\beta)$ there is a constant $C_{\alpha,\beta}$ with
\begin{equation}\label{Marcinkiewicz2}
|\xi|^{|\alpha|}|\eta|^{|\beta|}|\partial_{\xi}^{\alpha}
\partial_{\eta}^{\beta} \sigma (\xi,\eta)|\le
C_{\alpha,\beta}\big(\log(1+|\log\tfrac{|\xi|}{|\eta|}|
)\,\big)^{-\theta}\end{equation}
but
$W_{\sigma}$ is not of weak type $(p_1,p_2)$ for any $p_1,p_2>0.$

These examples indicate that the Marcinkiewicz-type conditions
(\ref{Marcinkiewicz}) need to be modified if they are to
imply boundedness for bilinear operators on $\mathbb R^{n} \times \mathbb R^{n}$.

In order to formulate some general results, let us introduce the following
notation. For $\sigma\in L_{\infty}$ we define
\begin{equation}\label{normH} \|\sigma\|_H=\sup_{1\le |\xi|\le
2}\sup_{1\le |\eta|\le 2}H((\sigma(2^j\xi,2^k\eta)_{j,k}).\end{equation}
If $\sigma$ is of class $C^N$ we define
\begin{equation}\label{normH2}
\|\sigma\|_{H}^{(N)}=\sum_{|\alpha|\le
N}\||\xi|^{|\alpha|}\partial_{\xi}^{\alpha}\sigma\|_H
+\sum_{|\beta|\le N}\||\eta|^{|\beta|}\partial_{\eta}^{\beta}\sigma\|_H.
\end{equation}
It will also be useful to define in this case
\begin{equation}\label{multiplierderiv}
\|\sigma\|_{\mathcal M_{p_1,p_2}}^{(N)}=
\sum_{|\alpha|\le
N}\||\xi|^{|\alpha|}\partial_{\xi}^{\alpha}\sigma\|_{\mathcal M_{p_1,p_2}}+
\sum_{|\beta|\le
N}\||\eta|^{|\beta|}\partial_{\eta}^{\beta}\sigma\|_{\mathcal M_{p_1,p_2}}.
\end{equation}

Now consider an arbitrary $L^{\infty}$ symbol $\sigma $ of class $C^{n+1}$.
Let
\begin{equation}\label{breakup}
\sigma_{jk}(\xi,\eta)=\sigma(\xi,\eta)\widehat\phi(2^{-j}\xi)
\widehat\phi(2^{-k}\eta).
\end{equation}
Set $\widehat{\zeta}(\xi)= \widehat{\phi_{-2}}(\xi)+
\widehat{\phi_{-3}}(\xi)+\widehat{\phi_{-4}}(\xi)$. Then $\widehat{\zeta}$ is
equal to $1$ on the annulus $1/16  \le |\xi|\le 1/4 $ and
vanishes off the annulus $1/32  \le |\xi|\le 1/2 $. Thus the function
$\widehat{\zeta}(\xi)\widehat{\zeta}(\eta)$ is supported in the unit cube
$[0,1]^{2n}$ and is equal to one on the support of
$$
(\xi,\eta) \to \sigma_{jk} (2^{j+3}\xi,2^{k+3}\eta)
$$
which is also contained in  $[0,1]^{2n}$. Inspired by \cite{CM2}, 
we expand the function above  in Fourier series on $[0,1]^{2n}$. We have
$$
\sigma_{jk} (2^{j+3}\xi,2^{k+3}\eta) =
\sum_{\nu\in \mathbb Z^n} \sum_{\rho\in \mathbb Z^n}
a_{jk}(\nu,\rho) e^{2\pi i (\langle \xi , \nu \rangle +
\langle \eta , \rho \rangle )} \widehat\zeta(\xi) \widehat\zeta(\eta),
$$
where   for $(\nu,\rho)\in\mathbb Z^n\times \mathbb Z^n $ we set
\begin{equation}\label{defmatrix}
a_{jk}(\nu,\rho) =  \int_{\mathbb R^n} \int_{\mathbb R^n}
\sigma(2^{j+3}t,2^{k+3}s)\widehat{\phi}(8 t)\widehat{\phi}(8s)
e^{-2\pi i( \langle t,\nu\rangle+
\langle s,\rho\rangle)} dt\,ds.
\end{equation}
We will denote by $A(\nu, \rho)$ the matrix with entries $a_{jk}(\nu,\rho)$.
Now setting
\begin{equation}\label{deftau}
\tau^{\nu,\rho}(\xi,\eta)=
\bigg(\sum_{j\in\mathbb Z}\sum_{k\in\mathbb
Z}a_{jk}(\nu,\rho)e^{ \frac{\pi i}{4} (2^{-j}\langle \xi,\nu\rangle+
2^{-k}\langle\eta,\rho\rangle)} \bigg)
 \widehat{\zeta}(2^{-j-3}\xi)  \widehat{\zeta }(2^{-k-3}\eta),
\end{equation}
we can write a symbol $\sigma$   of class $C^{n+1}$ as
\begin{equation}\label{pointwise}
\sigma(\xi,\eta)=\sum_{\nu\in\mathbb
Z^n}\sum_{\rho\in\mathbb Z^n}
\tau^{\nu,\rho}(\xi,\eta).
\end{equation}

In the next lemma we obtain some elementary estimates based on this
expansion.

\begin{Lem}\label{fourierestimates}
  Suppose $0<p_1,p_2<\infty$ and
$\frac1p_0=\frac1p_1+\frac1p_2.$  Then:\newline
(i) There is a constant $C=C(p_1,p_2)$ so
that for any $(\nu,\rho)$
$$
\|\tau^{\nu,\rho}\|_{\mathcal M_{p_1,p_2}}\le
C  (1+|\nu|+|\rho|)^{2m}H(A(\nu,\rho))
$$
where $m=[(n+1)/2].$\newline
(ii) There is a constant $C=C(N,p_1,p_2)$ such that if $\sigma$ is of class
$C^N,$ and $|\nu|+|\rho|>0$, then
$$
 H(A(\nu,\rho))\le C(1+|\nu|+|\rho|)^{-N}\|\sigma\|_H^{(N)},
$$
while
$$
H(A(0,0)) \le C\|\sigma\|_H.
$$
\newline (iii) If $p_0\ge 1$ and $\sigma$ is of class $C^N$ then there is a
constant $C=C(N,p_1,p_2)$ such that
$$
H(A(\nu,\rho))\le
C(1+|\nu|+|\rho)^{2m-N}\|\sigma\|_{\mathcal  M_{p_1,p_2}}^{(N)}.
$$
\end{Lem}

\begin{proof}
Observe  that $\widehat{\zeta}(2^{-j-3}\xi)= \widehat{\phi}(2^{-j-1}\xi) +
\widehat{\phi}(2^{-j }\xi) +\widehat{\phi}(2^{-j+1}\xi) $ and therefore
$\tau^{\nu,\rho}(\xi, \eta) $  is the sum of nine terms of the form
$$
\sum_{j\in\mathbb Z}\sum_{k\in\mathbb
Z}a_{j ,k }(\nu,\rho)  \big( e^{ \frac{\pi i}4 \langle 2^{-j}\xi,\nu\rangle}
 \widehat \phi(
2^{-j-r} \xi) \big)
\big(e^{ \frac{\pi i}4 \langle 2^{-k}\eta,\rho\rangle}\widehat  \phi
(2^{-k-s}  \eta)  \big)
$$
where $r,s\in \{-1,0,+1\}$.
We now use Lemma  \ref{factors},
Lemma \ref{translation2} (ii), and Lemma \ref{upper} in that order to obtain
$$
 \|\tau^{\nu,\rho} \|_{\mathcal M_{p_1,p_2}}\le
C(1+|\nu|)^m(1+|\rho|)^m H(A (\nu,\rho))
$$
where $m=[(n+1)/2].$  This proves (i).

For (ii) note that
if $|\alpha|,|\beta|\le N$   integration
by parts gives
 \begin{align}\label{fourier1}
 a_{jk}(\nu,\rho)&  \!=\!
\int_{\mathbb R^n}\!\int_{\mathbb R^n}\!\!\!
\partial_{\xi}^{\alpha} \big(\sigma(2^{j+3}\xi,2^{k+3}\eta)
\widehat{\phi}(8\xi)\widehat{\phi}(8\eta)\big)
 \frac{e^{-2\pi i( \langle \xi,\nu\rangle+
\langle\eta,\rho\rangle)}}{(-2\pi i \nu)^{\alpha}}
d\xi d\eta , \\  \label{fourier2}
 a_{jk}(\nu,\rho)& \!=\!
\int_{\mathbb R^n} \!\int_{\mathbb R^n}\!\!\!
\partial_{\eta}^{\beta}\big(\sigma(2^{j+3}\xi,2^{k+3}\eta)\widehat{\phi}(8\xi)
\widehat{\phi}(8\eta)\big)
\frac{e^{-2\pi i( \langle \xi,\nu\rangle+
\langle\eta,\rho\rangle)}}{(-2\pi i\rho)^{\beta}}
d\xi d\eta ,
\end{align}
provided   $\nu_1^{\alpha_1}\dots \nu_n^{\alpha_n}$ and $\rho^{\beta_1}\dots
\rho_n^{\beta_n}$ are nonzero.

Now using the fact that $H$ is a norm it is easy to see that by choosing
an appropriate $\alpha$ or $\beta$ for each pair $(\nu,\rho)\neq (0,0)$
one obtains the estimate
$$
H(A(\nu,\rho)) \le  C(N,p_1,p_2)(1+|\nu|+|\rho|)^{-N}\|\sigma\|_H^{(N)}.
$$
If $(\nu,\rho)=(0,0) $ the same estimate follows directly
from  (\ref{defmatrix}).

Finally we turn to (iii).
For fixed $\delta_j,\delta_k'$ with $\sup|\delta_j|,\sup|\delta_k'|\le 1$
let us define $\mu(\xi)=\sum_{j\in\mathbb Z}\delta_j\widehat{\phi_j}(\xi)$ and
$\upsilon(\eta)=\sum_{k\in\mathbb Z}\delta_k'\widehat{\phi_j}(\eta).$   Then it
follows from Lemma \ref{factors} that for any multi-indices
$\alpha,\alpha'$ we have
$$
\||\xi|^{|\alpha|+|\alpha'|}\partial_{\xi}^{\alpha}\mu(\xi)
\partial_{\xi}^{\alpha'} \sigma(\xi,\eta)\upsilon(\eta)\|_{\mathcal
M_{p_1,p_2}}\le C(\alpha,\alpha')\|\sigma\|^{(|\alpha'|)}_{\mathcal
M_{p_1,p_2}}$$
This implies that for fixed $N$ and any $\alpha$ with $|\alpha|=N$ we
have
\begin{equation}\label{a1}
\sup_{|\delta_j|\le 1}\sup_{|\delta_k'|\le 1}\||\xi|^{N}\sum_{j\in\mathbb
Z}\sum_{k\in\mathbb Z}
\delta_j\delta_k'\partial_{\xi}^{\alpha}\sigma_{jk}(\xi,\eta)
\|_{\mathcal
M_{p_1,p_2}}\le C(N)\|\sigma\|^{(N)}_{\mathcal
M_{p_1,p_2}}.
\end{equation}

We now use either  (\ref{defmatrix}) if $(\nu,\rho)=(0,0)$ or
we refer back to Proposition \ref{averaging}
(\ref{fourier1}) or (\ref{fourier2})
according to the values of $\nu$ or $\rho$,
when $(\nu,\rho) \neq (0,0)$.   For
example when $N=|\nu |\ge |\rho |$ and    the $l$th entry of $\nu$ has
maximal size $N$, then
\begin{align*}
&\|\sum_{j\in\mathbb Z}\sum_{k\in\mathbb
Z}a_{jk}(\nu,\rho)\widehat{\phi_j}(\xi)\widehat{\phi_k}(\eta)\|_{\mathcal
M_{p_1,p_2}}            \\
\le &C\sup_{|\delta_j|\le 1}\sup_{|\delta'_k|\le
1}\|\sum_{j\in \mathbb Z}\sum_{k\in\mathbb
Z}\delta_j\delta_k'2^{jN}\frac{\partial^{N}}{\partial\xi_l^N}
\sigma_{jk}(\xi,\eta) \frac{e^{-2\pi  i (\langle 2^{-j}
\xi,\nu\rangle+\langle 2^{-k} \eta,\rho\rangle)}}{(-2\pi i \nu_l)^N}\|_{\mathcal
M_{p_1,p_2}}.
\end{align*}
Now by Lemma \ref{factors} we can estimate the last expression side above by
$$
 C(1+|\nu|+|\rho|)^{2m-N}\sup_{|\delta_j|\le1}\sup_{|\delta_k'|\le 1}
 \|\sum_{j\in\mathbb Z}\sum_{k\in\mathbb Z}\delta_j\delta_k'|\xi|^N
\frac{\partial^{N}}{\partial\xi_l^N}
\sigma_{jk}(\xi,\eta)\|_{\mathcal M_{p_1,p_2}}.
$$
Using (\ref{a1}) we obtain (iii).
\end{proof}

Let us state the main result of this section.

\begin{Thm}\label{main1}
  Suppose $0<p_1,p_2<\infty$ and
$\frac1p_0=\frac1p_1+\frac1p_2.$  Let $N=2n+1$ if $p_0\ge 1$ and
$N=n+2+[\frac{n}{p_0}]$
if $p_0<1.$  Then for any $\sigma$   $C^N-$symbol such that
$\|\sigma\|_H^{(N)}<\infty$ we have $\|\sigma\|_{\mathcal
M_{p_1,p_2}}<\infty$.  Furthermore, there is a constant $C=C(p_1,p_2)$ so
that $\|\sigma\|_{\mathcal M_{p_1,p_2}}\le C\|\sigma\|_H^{(N)}.$
\end{Thm}

\begin{proof}   This follows directly from Lemma \ref{fourierestimates}
and (\ref{pointwise}).  Indeed, we have
$$ \|\tau^{\nu,\rho}\|_{\mathcal M_{p_1,p_2}}\le
C(1+|\nu|+|\rho|)^{2m-N}.$$  If $t=\min(p_0,1)$ we have
$$ \|\sigma\|_{\mathcal M_{p_1,p_2}}\le C (\sum_{\nu\in\mathbb
Z}\sum_{\rho\in\mathbb Z}
(1+|\nu|+|\rho|)^{(2m-N)t})^{\frac1t}\|\sigma\|_H^{(N)}.$$
Since $(N-2m)t>n$ this gives the result.\end{proof}

We next show that in a certain sense the preceding theorem is best
possible.

\begin{Thm}\label{bestposs} Suppose $1<p_1,p_2<\infty$ and
$\frac1p_0=\frac1p_1+\frac1p_2\le 1.$ Suppose
$\sigma$ is a
$C^{\infty}-$symbol.
 Then the following are equivalent:\newline
(i) $
\|\sigma\|_{\mathcal M_{p_1,p_2}}^{(N)}<\infty$ for every $N\ge
0.$\newline
(ii)$ \|\sigma\|_H^{(N)}<\infty$ for every $N\ge 0.$
\end{Thm}

\begin{proof}  Assume (i); then it follows from Lemma
\ref{fourierestimates} that for any $N>0$ we have an estimate
$H(A(\nu,\rho))\le C_N (1+|\nu|+|\rho|)^{-N}.$ Now it is clear from the
definition and from  Theorem \ref{matrixnorms} and Lemma
\ref{translation2} that we have an estimate
$$
 \||\xi|^{|\alpha|}\partial_{\xi}^{\alpha}\tau^{\nu,\rho}\|_H \le
C_{\alpha}(1+|\nu|)^{|\alpha|}H(A(\nu,\rho)).
$$  Hence we can deduce
easily that $$\||\xi|^{\alpha}\partial_{\xi}^{\alpha}\sigma\|_H
<\infty$$ for each multi-index $\alpha.$ Repeating the same reasoning
with the second variable $\eta$ gives
(ii).

Now assume (ii).  Then for any multi-index $\alpha$ one can see easily by
differentiation that for any pair of multi-indices $\alpha,\beta$ we have
that (ii) is satisfied by the symbols
$|\xi|^{|\alpha|}\partial_{\xi}^{\alpha}\sigma$ and
$|\eta|^{|\beta|}\partial_{\eta}^{\beta}\sigma$ in place of $\sigma.$
Applying Theorem \ref{main1} gives (i).
\end{proof}

Now let us recast Theorem \ref{main1} in terms of estimates on the symbol
$\sigma$ using the results of Section \ref{estimates}.

\begin{Thm}\label{main2} Suppose $0<p_1,p_2<\infty$ and
$\frac1p_0=\frac1p_1+\frac1p_2.$  Let $N=2n+1$ if $p_0\ge 1$ and
$N=n+2+[\frac{n}{p_0}]$
if $p_0<1.$  Suppose $\theta>1$
 Suppose
$\sigma$ is a
$C^N-$symbol such that for any pair of multi-indices $\alpha,\beta$ with
$0\le
|\alpha|\le N$ and $0\le |\beta|\le N$ there exist constants $C_{\alpha},
C_{\beta},$ with
\begin{align}\label{Marcinkiewicz3}
|\xi|^{|\alpha|} |\partial_{\xi}^{\alpha}\sigma(\xi,\eta)| &\le
C_{\alpha}
(\log(1+|\log\tfrac{|\xi|}{|\eta|}|))^{-\theta} \\
\label{Marcinkiewicz4}
|\eta|^{|\beta|} |\partial_{\eta}^{\beta}\sigma(\xi,\eta)| &\le
C_{\beta}
(\log(1+|\log\tfrac{|\xi|}{|\eta|}|))^{-\theta}.\end{align}
Then $\|\sigma\|_{\mathcal M_{p_1,p_2}}<\infty.$
\end{Thm}

\noindent {\bf Remark.}  We have already seen that in (\ref{Marcinkiewicz2}) that
this is false when $0<\theta<\frac12.$  However the arguments of Section
\ref{estimates} shows that we can improve  (\ref{Marcinkiewicz3}) and
(\ref{Marcinkiewicz4}) somewhat.  For example we can replace $
\big(\log(1+|\log\frac{|\xi|}{|\eta|}|)\big)^{-\theta} $ where $\theta>1$ by $
\big(  \log(1+|\log\frac{|\xi|}{|\eta|}|)\big)^{-1}\big(
\log(1+\log(1+|\log\frac{|\xi|}{|\eta|})\big)^{-\gamma}$  where $\gamma>1.$

\begin{proof}  This follows immediately from Theorem \ref{main1} and
Theorem \ref{lorentz} which yields the estimate
$$ H(A) \le C\sup_{j,k}\frac{|a_{jk}|}{w_{|j-k|+1}}$$ with $w_k=\log
(1+k)^{-\theta}.$\end{proof}

It is possible to ``mix and match'' the estimates in Section
\ref{estimates}: for example, in the following theorem we remove the
conditions for $|\alpha|,|\beta|=0$ but insist on a stronger condition
for
$|\alpha|=|\beta|=1$:

\begin{Thm}\label{main3}
Suppose $0<p_1,p_2<\infty$ and
$\frac1p_0=\frac1p_1+\frac1p_2.$  Let $N=2n+1$ if $p_0\ge 1$ and
$N=n+2+[\frac{n}{p_0}]$
if $p_0<1.$  Suppose $\theta>1$
 Suppose
$\sigma$ is a
$C^N-$symbol which satisfies conditions (\ref{Marcinkiewicz3}) and
(\ref{Marcinkiewicz4}) for $2\le |\alpha|,|\beta|\le N$ and if
$|\alpha|=|\beta|=1$
\begin{equation}\begin{align}\label{Marcinkiewicz5}
|\xi|^{|\alpha|} |\partial_{\xi}^{\alpha}\sigma(\xi,\eta)| &\le
C_{\alpha}
(1+|\log\tfrac{|\xi|}{|\eta|}|)^{-\theta} \\
\label{Marcinkiewicz6}
|\eta|^{|\beta|} |\partial_{\eta}^{\beta}\sigma(\xi,\eta)| &\le
C_{\beta}
(1+|\log\tfrac{|\xi|}{|\eta|}|)^{-\theta}.\end{align}\end{equation}
Then $\|\sigma\|_{\mathcal M_{p_1,p_2}}<\infty.$
\end{Thm}

\begin{proof} It is only necessary to show that $\|\sigma\|_H<\infty.$
Note first that Proposition \ref{bv} can be used to give the estimate
for any infinite matrix:
$$
H(A) \le C\big(\|A\|_{\infty}+\sup_j\sum_{k<j}|a_{j,k}-a_{j,k+1}|+
\sup_k\sum_{j<k}|a_{j,k}-a_{j+1,k}|\big).
$$
Now suppose $1\le |\xi|,|\eta|\le 2$.  Then if $k<j,$
$$ |\sigma(2^j\xi, 2^k\eta)-\sigma(2^j\xi,2^{k+1}\eta)|\le Ck^{-\theta}$$
by (\ref{Marcinkiewicz6}).  Combining with a similar estimate from
(\ref{Marcinkiewicz5}) gives the theorem.
\end{proof}

We conclude this section with a theorem of the type of Theorem
\ref{main1} for operators on $L_1.$

\begin{Thm}
\label{L1}  Suppose $N=2n+3$ and that $\sigma$ is a
$C^N$-symbol with $\|\sigma\|_H^{(N)}<\infty$; then
$W_{\sigma}:L_1\times L_1\to L_{\frac12,\infty}$ is bounded.
\end{Thm}

\begin{proof} Let $Q$ be the cube $\{x: \max_k|x_k|\le 1\}$ and
consider the bilinear operator
$W_{\sigma,Q}(f,g)=\chi_{Q}W_{\sigma}(f,g).$  We will show
that if $r<\frac12$ is such that $n+2+[\frac{n}{2r}]=N $,  then
$W_{\sigma,Q}:L_1(2Q)\times L_1(2Q)
\to L_r(Q)$ is bounded
and $\|W_{\sigma,Q}\|\le C\|\sigma\|_H^{(N)}$ where $C$ is a
constant depending only on dimension.

Suppose that $f,g\in\mathcal S$ are functions with support
contained in $2Q$ and such that
$\int f(x)\,dx=\int g(x)\,dx=0.$ Then $f,g\in H_{2r}$ with
$\|f\|_{H_{2r}}\le C\|f\|_{L_1}$ and $ \|g\|_{H_{2r}}\le C\|g\|_{L_1}.$  
Applying Theorem \ref{main1} we obtain that
\begin{equation}\label{req} 
\|W_{\sigma}(f,g)\|_{L_r} \le
C\|\sigma\|_H^{(N)}\|f\|_{L_1}\|g\|_{L_1}
\end{equation}  
where $C$ is an absolute constant.  It follows that $W_{\sigma}$ extends
unambiguously to any $f,g\in L_1(2Q)$ with $\int f(x)\,dx=\int
g(x)\,dx=0$ and (\ref{req}) holds.

 Next fix $\psi\in \mathcal S$ so that
$\int \psi(x)\,dx=1$ and $\psi$ has  support contained in $Q.$
Now for any $f,g\in L_1(3Q)$ let $f_0=f-(\int f(x)\,dx)\psi$ and
$g_0=g-(\int g(x)\,dx)\psi.$  Then   (\ref{req}) gives
$$
 \|W_{\sigma,Q}(f_0,g_0)\|_{L_r} \le C\|\sigma\|_H^{(N)}\|f\|_{L_1}\|g\|_{L_1}.
$$
We also note that $\|W_{\sigma,Q}(\psi,\psi)\|_{L_r} \le
C\|\sigma\|_H^{(N)}.$
Now consider the linear map $Tf= W_{\sigma}(f,\psi).$  Since $\psi\in
L_2$ we have that, if $\frac1s=\frac1{2r}+\frac12$,  $T:H_{2r}\to L_s$ is
bounded with norm controlled by $C\|\sigma\|_H^{(N)}$ (again using
Theorem \ref{main1}.)
Hence since $r<s,$
$$
 \|W_{\sigma,Q}(f_0,\psi)\|_{L_r} \le C\|\sigma\|_H^{(N)}\|f\|_{L_1}.
$$
Similarly
$$
 \|W_{\sigma,Q}(\psi,g_0)\|_{L_r} \le C\|\sigma\|_H^{(N)}\|g\|_{L_1}.
$$
Combining these estimates gives
\begin{equation}\label{start} \|W_{\sigma,Q}(f,g)\|_{L_r} \le
C\|\sigma\|_H^{(N)}\|f\|_{L_1}\|g\|_{L_1}.\end{equation}

We now use a Nikishin type argument as earlier in Lemma
\ref{Lem1.3}.
Suppose $(f_j)_{j=1}^J$ and $(g_j)_{j=1}^J$ satisfy $
\|f_j\|_{L_1},\|g_j\|_{L_1}\le 1$ and that $\sum_{j=1}^J|b_j|^{1/2}=1.$ Then if
$(\epsilon_j)_{j=1}^J$ and $(\epsilon'_j)_{j=1}^J$ are two independent
sequences of Bernoulli random variables we have
$$
 \big(\mathbb E(\|\sum_{j=1}^J\sum_{k=1}^J
\epsilon_j\epsilon'_k|b_j|^{\frac12}|b_k|^{\frac12}W_{\sigma,Q}
(f_j,g_k)\|_{L_r}^r) \big)^{\frac1r} \le C\|\sigma\|_H^{(N)}.
$$   
Again by using
the result of Bonami \cite {bonami}, we obtain an estimate
$$
 \|(\sum_{j=1}^J\sum_{k=1}^J |b_j||b_k|
|W_{\sigma,Q}(f_j,g_k)|^2)^{1/2}\|_{L_r} \le C\|\sigma\|_H^{(N)}.
$$
Extracting the diagonal gives
$$
 \big\|\max_{1\le j\le J}|b_j||W_{\sigma,Q}(f_j,g_j)|\big\|_{L_r} \le
C \|\sigma\|_H^{(N)}.
$$
We now use \cite{Pisier} as before.  There is a weight function $w\in
L_1(Q)$ with $w\ge 0$ a.e. and $\int w(x)\,dx=1$ so that for any $f,g\in
L_1(3Q)
$ with $\|f\|_{L_1},\|g\|_{L_1}\le 1$ and any measurable $E\subset Q$ we have
$$
\left (\int_E |W_{\sigma}(f,g)|^rdx\right)^{\frac1r} \le
C\|\sigma\|_H^{(N)} \left(\int_Ew(x)\,dx\right)^{\frac1r-2}.
$$
Now suppose $f,g$ are supported in $Q$ and $\lambda>0.$ Let $E=\{x\in Q:
|W_{\sigma}(f,g)|>\lambda.$  Then the above equation yields
\begin{equation}\label{mmmooo}
 \lambda |E|^{\frac1r} \le
C\|\sigma\|_H^{(N)}\left(\int_E w(x)\,dx\right)^{\frac1r-2}.
\end{equation}
On the other hand if we apply (\ref{mmmooo}) to $f_t(x)=f(x-t)$ where $t\in
Q$ and note that $W_{\sigma}(f_t,g)=(W_{\sigma}(f,g))_t$ we also obtain that
$$
 \lambda |E\cap (Q+t)|^{\frac1r} \le C\|\sigma\|_H^{(N)}
\left(\int_Ew(x-t)\,dx\right)^{\frac1r-2}.
$$
Raising to the power $(\frac1r-2)^{-1}$ and averaging gives:
$$
\lambda |E|^{\frac1r} \le C\|\sigma\|_H^{(N)} |E|^{\frac1r-2}.
$$
Thus $W_{\sigma,Q}$ maps $L_1(2Q)\times L_1(2Q)$ into
$L_{\frac12,\infty}(Q)$ with norm at most $C\|\sigma\|_H^{(N)}.$

Now let $\lambda>1$.  If we define
$\sigma_{\lambda}(\xi,\eta)=\sigma(\lambda^{-1}\xi,\lambda^{-1}\eta),$
then we have $\|\sigma_{\lambda}\|_H^{(N)}=\|\sigma_{\lambda}\|_H^{(N)}$
and we can apply this result to $\sigma_{\lambda}$. Notice that
$W_{\sigma_{\lambda}}(f,g)(x)=W_{\sigma}(f_{\lambda},g_{\lambda})(\lambda
x)$ where
$f_{\lambda}(x)=f(\lambda x)$ and $g_{\lambda}(x)= g(\lambda x).$
 This implies that
for any $\lambda>0$ we have the estimate
$$
 \|\chi_{\lambda Q}W_{\sigma}(f,g)\|_{L_{\frac12,\infty}}\le
C\|\sigma\|_H^{(N)}\|f\|_{L_1}\|g\|_{L_1}
$$ for $f,g$ supported in $\lambda Q$.
Letting $\lambda\to\infty$ gives the result. 
\end{proof}

\section{Discussion on paraproducts}\label{s-paraproducts}
 \setcounter{equation}0

Paraproducts are    bilinear operators of the type
$\sigma_A$ for some specific upper (or lower) triangular matrices $A$ of
zeros and ones. Paraproducts are important tools which have been used in several
occasions  in harmonic  analysis,
such as in the proof of the $T1$ theorem of David and Journ\'e  \cite{DJ}.
We define the  lower  and  upper  paraproducts  as the bilinear operators
$\Pi_L$ and $\Pi_U$ with symbols
$$
\tau_L(\xi,\eta)=\sum_{j\in\mathbb Z}
\sum_{k\le j-3}\widehat{\phi_j}(\xi)\widehat{\phi_k}(\eta)
$$
and
$$
\tau_U(\xi,\eta)=\sum_{k\in\mathbb
Z}\sum_{j\le k-3}\widehat{\phi_j}(\xi)\widehat{\phi_k}(\eta)
$$
respectively. It is easy to see that
$\|\tau_L\|_{\mathcal M_{p_1,p_2}},\|\tau_U\|_{\mathcal M_{p_1,p_2}}<\infty$ for all
$0<p_1,p_2<\infty$.  This can be deduced in several ways, e.g. from
Proposition \ref{subdiagonal} using Lemma \ref{translation2} or directly
from Theorem \ref{main1} and Proposition \ref{bv}.
We conclude that for all $0<p,q<\infty$
$\Pi_L$ maps $H_{p_1}\times H_{p_2}\to H_{p_0}$ when $1/p_1+1/p_2=1/p_0$
and $H_q=L_q$ when $1<q<\infty$.
We now turn to some endpoint cases regarding the paraproduct operator
$\Pi_L$.

\begin{Prop}
Let  $0<q<\infty$. Then the  paraproduct operator $\Pi_L$
is bounded on the following products of spaces.
\begin{enumerate}
\item[(1)] $  BMO\times H_q(\mathbb R^n)\to H_q(\mathbb R^n)$,
where $H_q=L_q$ when $1<q<\infty$.
\item[(2)] $ BMO\times H_1(\mathbb R^n)\to L_1(\mathbb R^n)$.
\item[(3)] $ BMO\times L_{\infty}(\mathbb R^n)\to BMO$.
\item[(4)] $  H_q(\mathbb R^n)\times L_{\infty}(\mathbb R^n)\to H_q(\mathbb
R^n)$, where $H_q=L_q$ when $1<q<\infty$.
\item[(5)] $ L_1(\mathbb R^n)\times L_{\infty}(\mathbb R^n)\to
L_{1,\infty}(\mathbb R^n)$.
\item[(6)] $  BMO\times  L_1(\mathbb R^n)\to L_{1,\infty}(\mathbb R^n)$.
\item[(7)] $  L_1(\mathbb R^n)\times  L_1(\mathbb R^n)\to L_{1/2,\infty}(\mathbb
R^n)$.
\end{enumerate}
\end{Prop}

\begin{proof}
Statement (1) is a classical result on paraproducts when $1<q<\infty$
and we refer the   reader to \cite{stein-new} p. 303 for a proof.
Note that for a fixed $f\in BMO$, the map
$g\to \Pi_L(f,g)$ is a Calder\'on-Zygmund singular integral.
The extension of (1) to $H_q$ for $q\le 1$, is   consequence of    the
that if a convolution type
singular integral operator maps $L_2\to L_2$ with bound
 a multiple of $\|f \|_{BMO}$, then it   also maps  $H_q$ into
itself with bound a multiple of this constant.   (2) follows from a similar
observation   while (3) is a dual statement to (2).
To prove (4) set $\widetilde S_jg= \sum_{k\le j-3} \widetilde \Delta_k g$.
We have that $\Pi_L(f,g)= \sum_{j \in \mathbb Z} \widetilde \Delta_j f
\widetilde S_jg$ and the Fourier transform of $\widetilde \Delta_j f
\widetilde S_jg$ is supported in the annulus $2^{j-2}\le |\xi| \le 2^{j+2}$.
It follows that
$$
\| \Pi_L(f,g) \|_{H_q}\le C \big\| \big(\sum_{j \in \mathbb Z}
| \widetilde \Delta_j f
\widetilde S_jg  |^2 \big)^{1/2}\big\|_{H_q} \le \|f\|_{H_q}\|Mg\|_{L_\infty},
$$
where $M$ is the Hardy-Littlewood maximal operator which is certainly
bounded on $L_\infty$.
To prove (5) we freeze   $g$   and look at the
linear operator $f\to \Pi_L(f,g)$ whose kernel is
$K(x,y)= \sum\limits_{j\in \mathbb Z} \phi_j(x-y) S_j(g)(x)$.
It is easy to see that
$$
|\nabla_{y} K(x,y)|\le  C\|g\|_{L_\infty} |x-y|^{-n-1}.
$$
This estimate together with the fact that the
linear operator $f\to \Pi_L(f,g)$ maps $L_2\to L_2$ gives that
$f\to \Pi_L(f,g)$ maps $L_1\to L_{1,\infty}$ using the
Calder\'on-Zygmund decomposition.  This proves (5).
 To obtain (6) we use  (1) (with $q=2$) and we apply to the
Calder\'on-Zygmund decomposition to the operator
$g\to \Pi_L(f,g)$ for fixed $f\in BMO$.
Finally (7) is a consequence of   Theorem \ref{L1}.
\end{proof}

\end{document}